# Elementary abelian subgroups of classical groups of type $A$

Meizheng Fu

ABSTRACT. Many open conjectures in the representation theory of finite groups can be studied by reducing them to related questions about quasi-simple groups. In such studies, $p$-radical subgroups typically play a critical role. To classify the $p$-radical subgroups of a finite group, we first classify the elementary abelian $p$-subgroups and find their local structure. To do this, we conduct the classification in a linear algebraic group $G$ and then transfer the results to the finite group of Lie type $G^F$. This approach was used by An, Dietrich and Litterick in their work on finite exceptional groups of Lie type. We now apply this approach to classical groups of type $A$ for the classification and local structure of the elementary abelian $p$-subgroups.

## 1. Introduction

One motivation for deep understanding of the elementary abelian subgroups of finite groups and algebraic groups of positive characteristic stems from profound questions arising from the study of finite simple groups. Griess [15] lists some of the reasons motivating his detailed study of the elementary abelian $p$-subgroups of algebraic groups. We mention just two examples. Using elementary abelian subgroups, Borel [9] studied torsion primes of connected compact Lie groups and established (using the classification of Lie groups) conditions for the existence of nontoral elementary abelian $p$-subgroups; Gorenstein and Lyons [13] investigated the self-centralizing elementary abelian subgroups of $K$-groups.

Elementary abelian subgroups also play a fundamental role in the representation theory of finite groups. Many open conjectures, such as the Alperin weight conjecture and Dade's invariant conjecture, can be studied by reducing to finite quasi-simple groups; see for instance [19]. Detailed knowledge of the local structure of $p$-radical subgroups facilitates this study. Here we refer to the normalizer and centralizer information as the local structure. A subgroup $R$ of a finite group $G$ is $p$-radical if $R = O_p(N_G(R))$, the largest normal $p$-subgroup of $N_G(R)$. A subgroup $H$ of $G$ is $p$-local, for some prime $p$, if there is some nontrivial $p$-subgroup $Q$ of $G$ with $H = N_G(Q)$. It is maximal $p$-local if it is $p$-local and maximal with respect to inclusion among all the $p$-local subgroups. If $M$ is $p$-local and maximal among all the proper $p$-local subgroups of $G$, then $M$ is a *maximal proper $p$-local* subgroup of $G$.

The classification of radical subgroups remains an open problem. One approach to classify $p$-radical subgroups of a finite group $G$ is via its $p$-local subgroups: every $p$-radical subgroup $R$ of a finite group $G$ with $O_p(G) \neq R$ is radical in some maximal proper $p$-local subgroup $M$ of $G$. Every such subgroup $M$ of $G$ can be realized as the normalizer of an elementary abelian $p$-subgroup. Thus, to classify the radical subgroups, we can first classify the elementary abelian subgroups. While this approach was applied successfully, in [5, 6] for example, to classify $p$-local and $p$-radical subgroups of some finite exceptional groups of Lie type, An, Dietrich and Litterick [4] adopted a different approach. There, the authors developed an algorithm which classifies the toral elementary abelian $p$-subgroups and their local structure in every simple algebraic group. (Hence we focus mainly on the nontoral ones.) Then they classified the nontoral elementary abelian $p$-subgroups and their local structure in a simple exceptional algebraic group $G$. Finally, they described the approach to descend to the finite simple exceptional group $G^F$ to determine the classification and local structure of its elementary abelian subgroups.

In this paper, we adopt the approach to carry out work in classical groups of type $A$. We only consider the elementary abelian $p$-subgroups of an algebraic group $G$ of type $A$ over $K$ with characteristic $\ell$ not $p$. By [14, Corollary 3.1.5], if $p = \ell$, then the normalizer of a $p$-radical subgroup of $G^F$ is a parabolic subgroup. We first recall the work of Andersen *et al.* [7] which provides the classification and local structure of nontoral elementary abelian $p$-subgroups of the Lie group $\mathrm{PGL}_n(\mathbb{C})$. We then classify the nontoral elementary abelian $p$-subgroups and identify the normalizer and centralizer of a representative of each class of the finite groups $\mathrm{PGL}_n(q)$ and $\mathrm{PGU}_n(q)$ where $p \nmid q$.

## 2. Notation

Throughout, unless stated otherwise, $G$ is a simple algebraic group defined over an algebraically closed field $K$ of characteristic $\ell \geq 0$. We view $G$ as a Zariski-closed subgroup of some general linear group over the algebraically closed field, except in Section 2.1 where we discuss briefly the scheme-theoretic background which is essential for transferring results among groups over fields with different characteristics. We fix a maximal torus $T$ of $G$. For a prime $p \neq \ell$, an elementary abelian $p$-subgroup of $G$ is *toral* if it has a conjugate contained in $T$ and *nontoral* otherwise. If $H$ is an algebraic group, then $H^\circ$ denotes the connected component of the identity element. Recall [18, Definition 21.3]: an endomorphism $F$ of $H$ is a *Steinberg endomorphism* if for some $m \geq 1$ the power $F^m : H \to H$ is the Frobenius morphism with respect to some $\mathbb{F}_{q^a}$-structure of $H$. An extension, split extension, and central product of a group $A$ by a group $B$ is written as $A.B$, $A : B$, and $A \circ B$, respectively. The direct product of $A$ and $B$ is denoted by $A \times B$, the cyclic group of order $m$ by $m$ or $C_m$, and a direct product of $n$ copies of $A$ by $A^n$. We denote a finite field with $q$ elements by GF($q$). For a $s \times s$ matrix $M$, we write $\Delta(M)_k$ to denote a $sk \times sk$ matrix with $k$ copies of $M$ embedded diagonally. As is defined in [4, Section 5.1], the *extended Weyl group* $\widetilde{W}$ of $G^F$ is an extension of $W$, the Weyl group of $G$, by an elementary abelian subgroup of order $\gcd(2, \ell-1)^r$, where $r$ is the rank of $G$, such that $\widetilde{W}$ maps onto $W$ under the natural projection map $N_G(T) \to W$. In particular $\widetilde{W} \leqslant G^F$, $\widetilde{W}$ normalizes $T^F$ and $\widetilde{W}T^F/T^F \cong W$.

## 3. Algebraic groups and finite groups

In this section, we prepare essential tools for our work. In Section 3.1, we introduce two results which allow us to transfer the classification and local structure from algebraic groups of characteristic 0 to those of positive characteristic. In Section 3.2, we record results which allow us to obtain the splittings in the finite group of a class in the algebraic group and the local structure of the splittings.

### 3.1. The defining characteristic of the algebraic group

The following propositions allow us to classify the elementary abelian $p$-subgroups of an algebraic group and obtain their local structure independently of the characteristic of the underlying field of the algebraic group so long as it is not $p$.

**Proposition 3.1** ([16, Theorem 1.22]). *Let $\ell$ be a prime number and $G(\cdot)$ a semisimple split group scheme. If $M$ is a finite group of order prime to $\ell$, then for algebraically closed fields $K_0$ and $K_\ell$ of characteristic $0$ and $\ell$, respectively, the sets $\mathrm{Hom}(M, G(K_0))/G(K_0)$ and $\mathrm{Hom}(M, G(K_\ell))/G(K_\ell)$ have the same cardinality. Furthermore, the process of reduction modulo $\ell$ induces a bijection.*

The process of reduction modulo $\ell$' entails showing that every finite subgroup of $G(K_0)$ has a conjugate contained in $G(R)$, where $R$ is a suitable ring of Witt vectors that maps onto an algebraically closed subfield of $K_\ell$, and then considering the induced map $G(R) \to G(K_\ell)$. We refer to [16, Appendix A] for more details.

**Proposition 3.2** ([4, Lemma 3.3]). *Let $K_0$ and $K_\ell$ be algebraically closed fields of characteristic $0$ and $\ell$, respectively and let $p \neq \ell$ be a prime. Let $E$ and $E'$ respectively be elementary abelian $p$-subgroups of $G(K_0)$ and $G(K_\ell)$ which correspond under the bijection of Proposition 3.1. The process of reduction modulo $\ell$ allows us to identify the root data of the reductive groups $C_{G(K_0)}(E)^\circ$ and $C_{G(K_\ell)}(E')^\circ$, and also induces the following isomorphisms:*

$$N_{G(K_0)}(E)/C_{G(K_0)}(E) \cong N_{G(K_\ell)}(E')/C_{G(K_\ell)}(E')$$
$$N_{G(K_0)}(E)/C_{G(K_0)}(E)^\circ \cong N_{G(K_\ell)}(E')/C_{G(K_\ell)}(E')^\circ$$
$$C_{G(K_0)}(E)/C_{G(K_0)}(E)^\circ \cong C_{G(K_\ell)}(E')/C_{G(K_\ell)}(E')^\circ.$$

### 3.2. From algebraic groups to finite groups

Once it is clear how to classify the elementary abelian subgroups of an algebraic group and find their local structure, we can tackle the finite group case. Let $G$ be a connected reductive algebraic group over an algebraically closed field $K$ of characteristic $\ell$, with maximal torus $T$. We suppose that a Steinberg endomorphism $F$ of $G$ is given and that $T$ is $F$-stable.

**Definition 3.3.** *If $H$ is a subgroup of $G$, then the $F$-class of a coset $aN_G(H)^\circ$ in $N_G(H)/N_G(H)^\circ$ is $\{F(g)ag^{-1}N_G(H)^\circ : g \in N_G(H)\}$.*

If $H$ is finite, then $N_G(H)^\circ = C_G(H)^\circ$. We use the following proposition to obtain information on the class distributions of the elementary abelian subgroups when we descend from algebraic groups to finite groups.

**Proposition 3.4** ([4, Proposition 4.1])**.** *Suppose that $G$ is a connected reductive algebraic group over an algebraically closed field of characteristic $\ell$, that a Steinberg endomorphism $F$ of $G$ is given, and that $A$ is a finite subgroup of $G$ where $A$ and $F(A)$, the image of $A$ under $F$, are conjugate in $G$.*
  (i) *There exists a $G$-conjugate of $A$ which is $F$-stable.*
  (ii) *Suppose $A$ has a conjugate in $G^F$. Replacing $A$ by this conjugate, there is a bijection between $G^F$-classes of subgroups of $G^F$ which are $G$-conjugate to $A$ and $F$-classes in $N_G(A)/C_G(A)^\circ$ contained in $C_G(A)/C_G(A)^\circ$: the $F$-class of $w \in C_G(A)/C_G(A)^\circ$ corresponds to the $G^F$-class of subgroups with representative $A_w = gAg^{-1}$, where $g \in G$ is chosen with $g^{-1}F(g)C_G(A)^\circ = w$.*

The following result tells when $A$ has a $G$-conjugate in $G^F$.

**Corollary 3.5** ([4, Corollary 4.2])**.** *Let $A$ be a finite subgroup of $G$.*
  (i) *There is a $G$-conjugate of $A$ in $G^F$ if and only if $A$ and $F(A)$ are $G$-conjugate and there is an $F$-stable $G$-conjugate $B$ of $A$ such that the restriction $F : B \to B$ is induced by some $w \in N_G(B)/C_G(B)^\circ$.*
  (ii) *If $|N_G(A)/C_G(A)| = |\mathrm{Aut}(A)|$, then there is a $G$-conjugate of $A$ in $G^F$ if and only if $A$ and $F(A)$ are $G$-conjugate.*

Once the conjugacy classes of the elementary abelian subgroups of the finite group are decided, we apply the following propositions to deduce their local structure.

**Proposition 3.6** ([4, Proposition 4.3])**.** *Let $A$ be a finite subgroup of $G^F$. For $w \in C_G(A)/C_G(A)^\circ$ let $A_w \leqslant G^F$ be the $G$-conjugate of $A$ as defined in Proposition 3.4. If $\dot{w} \in C_G(A)$ is a preimage of $w$, then*
$$(C_G(A_w)^\circ)^F \cong (C_G(A)^\circ)^{\dot{w}F}$$
*where the map $\dot{w}F$ sends $x$ to $\dot{w}F(x)\dot{w}^{-1}$, and $(C_G(A)^\circ)^{\dot{w}F}$ are the fixed points of $\dot{w}F$ in $C_G(A)^\circ$. Furthermore, $(C_G(A)^\circ)^{\dot{w}F}$ is independent of the choice of $\dot{w}$. If $w$ acts as an inner automorphism of $C_G(A)^\circ$, then $(C_G(A_w)^\circ)^F \cong (C_G(A)^\circ)^F$.*

**Proposition 3.7** ([4, Proposition 4.4])**.** *Let $A$ be a finite subgroup of $G^F$. For $w \in C_G(A)/C_G(A)^\circ$ let $A_w \leqslant G^F$ be the $G$-conjugate of $A$ as defined in Proposition 3.4.*

$$N_G(A_w)^F/(C_G(A_w)^\circ)^F \cong (N_G(A)/C_G(A)^\circ)^{wF}$$
$$C_G(A_w)^F/(C_G(A_w)^\circ)^F \cong (C_G(A)/C_G(A)^\circ)^{wF}.$$

The following corollary is useful when $F$ acts trivially on $N_G(A)/C_G(A)^\circ$.

**Corollary 3.8** ([4, Corollary 4.7])**.** *Let $A$ be a finite subgroup of $G^F$. If the Steinberg endomorphism $F$ acts trivially on $N_G(A)/C_G(A)^\circ$, then there is a bijection between the $N_G(A)/C_G(A)^\circ$-classes in $C_G(A)/C_G(A)^\circ$ and the $G^F$-classes of subgroups of $G^F$ that are $G$-conjugate to $A$; this bijection maps $w \in C_G(A)/C_G(A)^\circ$ to $A_w \leqslant G^F$, as defined in the correspondence in Proposition 3.4. Moreover*

$$C_{G^F}(A_w) = (C_G(A_w))^F = (C_G(A)^\circ)^{\dot{w}F}.C_{C_G(A)/C_G(A)^\circ}(w)$$
$$N_{G^F}(A_w) = (N_G(A_w))^F = (C_G(A)^\circ)^{\dot{w}F}.C_{N_G(A)/C_G(A)^\circ}(w).$$

## 4. Elementary abelian $p$-subgroups of classical groups of type $A$

Let $K$ be an algebraically closed field of characteristic $\ell$ not $p$. Now we investigate the elementary abelian $p$-subgroups and their local structure in classical groups of type $A$. Since $\text{SL}_n(K)$ has no nontoral elementary abelian $p$-subgroups (see [20]), we work in the algebraic group $\text{PGL}_n(K)$ and then descend to the finite group.

In the following, we first study elementary abelian 2-subgroups and then elementary abelian $p$-subgroups where $p$ is an odd prime. Section 4.1 develops required theory. Section 4.3 is about nontoral elementary abelian 2-subgroups of $\text{PGL}_n(q)$ where $q \equiv 1 \pmod 4$. Section 4.4 deals with the case where $q$ is a power of a prime $\ell$ and $\ell \equiv 3 \pmod 4$. In Section 4.5, we consider the elementary abelian 2-subgroups of $\text{PGU}_n(q)$. Sections 4.6 and 4.7 examines the elementary abelian $p$-subgroups where $p$ is odd.

We use $[A]$ to denote the image of $A \in \text{GL}_n(\mathbb{C})$ in $\text{PGL}_n(\mathbb{C})$. A bar over a coset representative denotes the coset and $\lfloor n \rfloor$ is the floor of $n$.

### 4.1. Nontoral elementary abelian $p$-subgroups of $\text{PGL}_n(K)$

In this section, we introduce the group $\bar{\Gamma}_r$ first. Then we establish a fundamental result critical to our investigation in groups of type $A$.

**Proposition 4.1** ([4, Lemma 3.2])**.** *Let $G$ be a reductive complex linear algebraic group with a maximal compact subgroup $H$, and let $E$ be an elementary abelian subgroup of $H$.*

(i) *The inclusion $H \hookrightarrow G$ induces a bijection between the conjugacy classes of the elementary abelian subgroups of $H$ and the conjugacy classes of the elementary abelian subgroups of $G$.*

(ii) *The group $C_G(E)^\circ \cap H$ is a maximal compact subgroup of $C_G(E)^\circ$ and has the same root datum as $C_G(E)^\circ$.*

(iii) $N_H(E) = (H \cap (N_G(E))^\circ).(N_G(E)/N_G(E)^\circ)$;
$C_H(E) = (H \cap (C_G(E))^\circ).(C_G(E)/C_G(E)^\circ)$.

This proposition allows us to view the classification of elementary abelian subgroups of $\mathrm{PGL}_n(\mathbb{C})$ in $\mathrm{PU}_n(\mathbb{C})$, the maximal compact subgroup.

Next we introduce an important group which plays a critical role in the classification of nontoral elementary abelian subgroups.

Let $S_{p^r}$ be the symmetric group on $p^r$ points. For $s = 0, 1, ..., r-1$, let $\sigma_s \in S_{p^r}$ be defined by

$$\sigma_s(i) = \begin{cases} i + p^s & \text{if } i \equiv 1, ..., (p-1)p^s \mod p^{s+1} \\ i - (p-1)p^s & \text{if } i \equiv (p-1)p^s + 1, ..., p^{s+1} \mod p^{s+1}. \end{cases}$$

Let $\beta = e^{2\pi i/p}$ be a primitive $p$th root of unity. Define matrices $A_0, ..., A_{r-1}, B_0, ..., B_{r-1}$ in the compact unitary Lie group $\mathrm{U}_{p^r}(\mathbb{C})$ of degree $p^r$ as

$$(A_s)_{ij} = \begin{cases} \beta^{[(i-1)/p^s]} & \text{if } i = j \\ 0 & \text{if } i \neq j \end{cases}$$

and

$$(B_s)_{ij} = \begin{cases} 1 & \text{if } \sigma_s(i) = j \\ 0 & \text{if } \sigma_s(i) \neq j \end{cases}$$

for $s = 0, 1, ..., r-1$. Therefore the $A_s$'s are diagonal and the $B_s$'s are permutation matrices. These matrices satisfy the following relations:

$$[A_t, A_s] = I = [B_t, B_s] = [B_t, A_s] \, (t \neq s), \, [B_t, A_t] = \beta \cdot I.$$

The matrices $A_0, ..., A_{r-1}$ and $B_0, ..., B_{r-1}$ generate elementary abelian $p$-subgroups $\mathbf{A}_r$ and $\mathbf{B}_r$ of rank $r$, respectively.

**Lemma 4.2** ([8, Corollary 6.3]). *The subgroups $\mathbf{A}_r$ and $\mathbf{B}_r$ are conjugate in $\mathrm{U}_{p^r}(\mathbb{C})$.*

By this lemma and Proposition 4.1, the images $\bar{\mathbf{A}}_r$ and $\bar{\mathbf{B}}_r$ of $\mathbf{A}_r$ and $\mathbf{B}_r$ are conjugate in $\mathrm{PGL}_n(\mathbb{C})$ where $n = p^r$.

**Definition 4.3.** *For a prime $p$ and positive integer $r$,*
$$\Gamma_r = \langle u \cdot I, A_s, B_s \mid u \in \mathbb{C}^\times, 0 \leq s \leq r-1 \rangle \leqslant \mathrm{U}_{p^r}(\mathbb{C}).$$

We let $\bar{\Gamma}_r = \langle \bar{A}_s, \bar{B}_s \mid 0 \leq s \leq r-1 \rangle$ where $\bar{A}_s$ and $\bar{B}_s$ are the images of $A_s$ and $B_s$ in $\mathrm{PU}_{p^r}(\mathbb{C})$. Clearly, $\bar{\Gamma}_r \cong p^{2r}$.

When $n = p^r k$, we take $k$ copies of each generator of $\bar{\Gamma}_r$ to form $2r$ block diagonal matrices in $\mathrm{PU}_n(\mathbb{C})$. To simplify notation, we also denote by $\bar{\Gamma}_r$ the subgroup of $\mathrm{PU}_n(\mathbb{C})$ generated by these matrices. Next we record the key theorem needed to determine the nontoral elementary abelian subgroups of $\mathrm{PGL}_n(\mathbb{C})$ and their local structure.

**Theorem 4.4** ([7, Theorem 8.5]). *Suppose $E$ is a nontoral elementary abelian $p$-subgroup of $\mathrm{PGL}_n(\mathbb{C})$ for an arbitrary prime $p$.*
  (i) *Up to conjugacy, $E$ can be written as $\bar{\Gamma}_r \times \bar{A}$, for some $r \geq 1$ with $n = p^r k$ and some abelian subgroup $A$ of $C_{\mathrm{GL}_n(\mathbb{C})}(\Gamma_r) \cong \mathrm{GL}_k(\mathbb{C})$.*
  (ii) *For a given $r$, the conjugacy classes of such subgroups $E$ are in one-to-one correspondence with the conjugacy classes of toral elementary abelian $p$-subgroups $\bar{A}$ of $\mathrm{PGL}_k(\mathbb{C}) \cong C_{\mathrm{PGL}_n(\mathbb{C})}(\bar{\Gamma}_r)^\circ$ (allowing the trivial subgroup), and the centralizer of $E$ is given by $C_{\mathrm{PGL}_n(\mathbb{C})}(E) \cong \bar{\Gamma}_r \times C_{\mathrm{PGL}_k(\mathbb{C})}(\bar{A})$.*
  (iii) *The Weyl group*
$$\begin{aligned}
W_{\mathrm{PGL}_n(\mathbb{C})}(E) &:= N_{\mathrm{PGL}_n(\mathbb{C})}(E)/C_{\mathrm{PGL}_n(\mathbb{C})}(E) \\
&= \begin{pmatrix} \mathrm{Sp}_{2r}(p) & 0 \\ p^{\mathrm{rk}\bar{A} \times 2r} & W_{\mathrm{PGL}_k(\mathbb{C})}(\bar{A}) \end{pmatrix} \\
&= p^{\mathrm{rk}\bar{A} \times 2r} : (\mathrm{Sp}_{2r}(p) \times W_{\mathrm{PGL}_k(\mathbb{C})}(\bar{A})).
\end{aligned}$$
  (iv) $\alpha \in \mathrm{Sp}_{2r}(p)$ *acts up to conjugacy as $\alpha \times 1$ on $C_{\mathrm{PGL}_n(\mathbb{C})}(E)$.*

The matrix form of $W_{\mathrm{PGL}_n(\mathbb{C})}(E)$ is used when we consider the group action. By Schur's Lemma, to embed $e \in \bar{A}$ in $\mathrm{PGL}_n(\mathbb{C})$, we replace each entry $a_{ij}$ of $e$ by $a_{ij} \cdot I_{p^r}$.

## 4.2. Maximal nontoral 2-subgroups and their local structure

In this section, we study the maximal nontoral elementary abelian 2-subgroups and determine their local structure. In Section 4.2.1, we obtain the structure and local structure of such groups using Theorem 4.4. In Section 4.2.2, we determine the $F$-action on $N/C^\circ$. This lays the foundation for the study in Section 4.3 of the nonmaximal nontoral elementary abelian 2-subgroups.

### 4.2.1. *The maximal nontoral 2-subgroups and their local structure*

Let $\Lambda$ be a maximal nontoral elementary abelian 2-subgroup of $G = \mathrm{PGL}_n(K)$ where $K$ is an algebraically closed field of characteristic $\ell$ where $\ell \equiv 1 \bmod 4$. Let $F$ be a Steinberg endomorphism such that $\Lambda \leq G^F = \mathrm{PGL}_n(q)$ where $q$ is a power of $\ell$. Denote $N_{\mathrm{PGL}_n(K)}(\Lambda)$ by $N$ and $C_{\mathrm{PGL}_n(K)}(\Lambda)$ by $C$.

Let $E \leq \mathrm{PGL}_n(\mathbb{C})$ correspond to $\Lambda$ under the bijection of Proposition 3.1. By Theorem 4.4, if $n = 2^s \times t$ where $\gcd(2, t) = 1$ and $s \geq 1$, then there are $s$ conjugacy classes of maximal nontoral elementary abelian 2-subgroups $E$ of $\mathrm{PGL}_n(\mathbb{C})$ which take the shape $\bar{\Gamma}_i \times 2^m$ where $1 \leq i \leq s$, $m = 2^{s-i}t - 1$ and $2^m$ is a representative of the conjugacy class of the maximal toral elementary abelian 2-subgroups of $\mathrm{PGL}_{m+1}(\mathbb{C})$. Every nontoral elementary abelian 2-subgroup is contained in a maximal one.

**(1)** When $n = 2^r \times k$ where $r \geq 1$, $k > 1$ and $n$ is not a power of 2, the maximal nontoral elementary abelian 2-subgroups of $\mathrm{PGL}_n(\mathbb{C})$ are $E = \bar{\Gamma}_r \times 2^m$ where $m = k - 1$. By Theorem 4.4, $C_{\mathrm{PGL}_n(\mathbb{C})}(E) \cong \bar{\Gamma}_r \times T_m$ and $W_{\mathrm{PGL}_n(\mathbb{C})}(E) = 2^{m \times 2r} : (\mathrm{Sp}_{2r}(2) \times S_{m+1})$.

**(2)** When $n = 2^s$ where $s \geq 2$, the representatives of the conjugacy classes of the maximal nontoral elementary abelian 2-subgroups of $\mathrm{PGL}_n(\mathbb{C})$ take three shapes. They are listed in Table 1.

TABLE 1. Maximal nontoral elementary abelian 2-subgroups of $U = \mathrm{PGL}_{2^s}(\mathbb{C})$

| Maximal nontoral $E$ | $C_U(E)$ | $C_U(E)/C_U(E)^\circ$ | $N_U(E)/C_U(E)$ |
|---|---|---|---|
| $\bar{\Gamma}_r \times 2^{2^{s-r}-1}$, $1 \leq r \leq s-2$ | $\bar{\Gamma}_r \times T_{2^{s-r}-1}$ | $\bar{\Gamma}_r$ | $2^a : (\mathrm{Sp}_{2r}(2) \times S_{2^{s-r}})$ |
| $\bar{\Gamma}_{s-1} \times 2$ | $\bar{\Gamma}_{s-1} \times (T_1.2)$ | $\bar{\Gamma}_{s-1} \times 2$ | $2^b : \mathrm{Sp}_{2(s-1)}(2)$ |
| $\bar{\Gamma}_s$ | $\bar{\Gamma}_s$ | $\bar{\Gamma}_s$ | $\mathrm{Sp}_{2s}(2)$ |

The 2 in $(T_1.2)$ equals $\left\langle \left[ \begin{pmatrix} 0 & I_{2^{s-1}} \\ I_{2^{s-1}} & 0 \end{pmatrix} \right] \right\rangle$, $a = (2^{s-r} - 1) \times 2r$ and $b = 1 \times 2(s-1)$.

4.2.2. *The $F$-action on the normalizer quotients of these subgroups.*

First, we record a lemma.

**Lemma 4.5** ([4, Lemma 5.2]). *Let $G$ be a connected reductive simple algebraic group and let $L$ be an elementary abelian $p$-subgroup of $G$. Let $T$ be a maximal torus of $G$, let $\widetilde{W}$ be an extended Weyl subgroup of $G$ and let $F$ be a Steinberg endomorphism of $G$. If $L \leq T \cap G^F$ and $g \in N_G(L) \setminus C_G(L)^\circ$, then $g = vc$ for some $v \in \widetilde{W}$ and $c \in C_G(L)^\circ$. In particular,*
$$N_G(L)/C_G(L)^\circ = N_{G^F}(L)/(C_G(L)^\circ)^F$$
*and*
$$C_G(L)/C_G(L)^\circ = C_{G^F}(L)/(C_G(L)^\circ)^F.$$
*Hence $F$ induces the identity map on $N_G(L)/C_G(L)^\circ$.*

Proposition 3.4 states that there is a bijection between the $\mathrm{PGL}_n(q)$-classes of subgroups of $\mathrm{PGL}_n(q)$ which are $\mathrm{PGL}_n(K)$-conjugate to $\Lambda$ and the $F$-classes of $N/C^\circ$ contained in $C/C^\circ$. So we focus on the $F$-classes of $N/C^\circ$ contained in $C/C^\circ$.

By Lemma 4.5, when $\Lambda$ is toral, the Steinberg endomorphism $F$ acts trivially on $N/C^\circ$. So the $F$-classes are just the conjugacy classes. By [10, Section 8.2] and [17, Proposition 4.6.6], for $r \geq 1$, $\bar{\Gamma}_r.\mathrm{Sp}_{2r}(2) \leq N_G(\Lambda)$ from **(1)** and **(2)** is a subgroup of $\mathrm{PGL}_{2^r}(q)$ for $q$ a power of a prime $\ell$ where $\ell \equiv 1 \bmod 4$. Thus, $\bar{\Gamma}_r.\mathrm{Sp}_{2r}(2)$ is centralized by the Steinberg endomorphism $F$.

Let $U = 2^{m \times 2r} : (\mathrm{Id} \times \mathrm{Id}) \leq W_{\mathrm{PGL}_n(\mathbb{C})}(E)$ from **(1)**. We denote the group $2^{m \times 2r}$ by $R$. Let $A \in \mathrm{Sp}_{2r}(2)$, $B \in S_{m+1}$ and $u \in R$. Now

$$\begin{pmatrix} A & 0 \\ 0 & B \end{pmatrix} \begin{pmatrix} \mathrm{Id} & 0 \\ u & \mathrm{Id} \end{pmatrix} \begin{pmatrix} A^{-1} & 0 \\ 0 & B^{-1} \end{pmatrix} = \begin{pmatrix} \mathrm{Id} & 0 \\ BuA^{-1} & \mathrm{Id} \end{pmatrix}$$

implies that the generators of $R$ are in one orbit under the action of $\mathrm{Sp}_{2r}(2) \times S_{m+1}$. Observe that $R$ has a generator

$$\left[ \begin{pmatrix} \Delta(\begin{pmatrix} -1 & \\ & 1 \end{pmatrix}))_{2^r-1} & \\ & \Delta(I_{2^r})_m \end{pmatrix} \right] C_{\mathrm{PGL}_n(\mathbb{C})}(E).$$

Since all of the above nonzero entries are 1 or $-1$, when $n$ is not a power of 2 and $\Lambda$ is a maximal nontoral elementary abelian 2-subgroup of $G$, the Steinberg endomorphism $F$ centralizes $N/C^\circ$. This conclusion also holds for **(2)**.

Our above construction provides $2mr$ elements $x_{ij}, y_{ij} \in \mathrm{PGL}_n(\mathbb{C})$ for $1 \leq i \leq m$ and $1 \leq j \leq r$ where

$$x_{ij} = \left[ \begin{pmatrix} \Delta(I_{2^r})_{i-1} & & \\ & \Delta(\begin{pmatrix} & -I_{2^{j-1}} \\ I_{2^{j-1}} & \end{pmatrix}))_{2^{r-j}} & \\ & & \Delta(I_{2^r})_{m-(i-1)} \end{pmatrix} \right]$$

$$y_{ij} = \left[ \begin{pmatrix} \Delta(I_{2^r})_{i-1} & & \\ & \Delta(\begin{pmatrix} -I_{2^{j-1}} & \\ & I_{2^{j-1}} \end{pmatrix}))_{2^{r-j}} & \\ & & \Delta(I_{2^r})_{m-(i-1)} \end{pmatrix} \right].$$

Recall $\bar{\Gamma}_r = \langle \bar{A}_s, \bar{B}_s \mid 0 \leq s \leq r-1 \rangle$ from Section 4.1. Let $J \cong 2^m$ from **(1)** be the representative of the conjugacy class of the maximal toral elementary abelian 2-subgroups of rank $m$ of $\mathrm{PGL}_k(\mathbb{C})$. The generators of $J$ are $g_i = [\mathrm{M}_i]$ where $i = 1, ..., m$ and $\mathrm{M}_i$ is a $k \times k$ diagonal matrix with $-1$ in the $i$-th position and 1 in the remaining positions. Hence the maximal nontoral subgroup $E = \bar{\Gamma}_r \times J$ of $\mathrm{PGL}_n(\mathbb{C})$ has generators $\bar{A}_s$, $\bar{B}_s$ and $g_i \otimes I_{2^r}$ for $0 \leq s \leq r-1$ and $i = 1, ..., m$.

Observe

$$\begin{pmatrix} & -1 \\ 1 & \end{pmatrix} \begin{pmatrix} 1 & \\ & -1 \end{pmatrix} \begin{pmatrix} & 1 \\ -1 & \end{pmatrix} = \begin{pmatrix} -1 & \\ & 1 \end{pmatrix},$$

$$\begin{pmatrix} & -1 \\ 1 & \end{pmatrix} \begin{pmatrix} & 1 \\ 1 & \end{pmatrix} \begin{pmatrix} & 1 \\ -1 & \end{pmatrix} = \begin{pmatrix} & -1 \\ -1 & \end{pmatrix}$$

and
$$\begin{pmatrix} -1 & \\ & 1 \end{pmatrix} \begin{pmatrix} & 1 \\ 1 & \end{pmatrix} \begin{pmatrix} -1 & \\ & 1 \end{pmatrix} = \begin{pmatrix} & -1 \\ -1 & \end{pmatrix}.$$

Hence

$$\begin{aligned}
x_{ij}\bar{A}_s x_{ij}^{-1} &= \begin{cases} \bar{A}_s & s+1 \neq j \\ \bar{A}_s(g_i \otimes I_{2^r}) & s+1 = j \end{cases} \\
x_{ij}\bar{B}_s x_{ij}^{-1} &= \begin{cases} \bar{B}_s & s+1 \neq j \\ \bar{B}_s(g_i \otimes I_{2^r}) & s+1 = j \end{cases} \\
y_{ij}\bar{A}_s y_{ij}^{-1} &= \bar{A}_s \text{ for all } s \\
y_{ij}\bar{B}_s y_{ij}^{-1} &= \begin{cases} \bar{B}_s & s+1 \neq j \\ \bar{B}_s(g_i \otimes I_{2^r}) & s+1 = j \end{cases}
\end{aligned}$$

$x_{ij}$ and $y_{ij}$ fix all the $g_i \otimes I_{2^r}$ for $i = 1,...,m$.

Let
$$N_1 = \langle x_{ij}, y_{ij} \mid 1 \leq i \leq m \text{ and } 1 \leq j \leq r \rangle.$$

If $n \in N_1$ is nontrivial, then $n \in N_{\mathrm{PGL}_n(\mathbb{C})}(E) \setminus C_{\mathrm{PGL}_n(\mathbb{C})}(E)$. Due to the action of $N_1$ on $E$ and the structure of $W_{\mathrm{PGL}_n(\mathbb{C})}(E)$, we can identify $N_1 C_{\mathrm{PGL}_n(\mathbb{C})}(E)/C_{\mathrm{PGL}_n(\mathbb{C})}(E)$ with $R$, the subgroup of $W_{\mathrm{PGL}_n(\mathbb{C})}(E)$. Observe that $C_G(\Lambda)/C_G(\Lambda)^\circ \cong \bar{\Gamma}_r$ is centralized by $N_1 T_m / T_m$. This construction leads to the same conclusion for the first and third cases of (2).

### 4.3. Nonmaximal nontoral 2-subgroups of $\mathrm{PGL}_n(q)$ where $q \equiv 1 \bmod 4$

Here we study the nonmaximal nontoral elementary abelian 2-subgroups $E = \bar{\Gamma}_r \times \bar{A}$. In Section 4.3.1 we first investigate the $F$-action on the normalizer quotients of $E$ in the algebraic groups. In Section 4.3.2 we study the structure and local structure of $\bar{A}$ where $\bar{A}$ has a disconnected centralizer. Finally, in Section 4.3.3 we descend to the finite groups and discuss separately the cases where $\bar{A}$ has a connected and disconnected centralizer.

4.3.1. *The F-action on the normalizer quotients of nontoral 2-subgroups.*

Let $K$ be an algebraically closed field of characteristic $\ell$ where $\ell \equiv 1 \bmod 4$. Let $F$ be a Steinberg endomorphism such that $\mathrm{PGL}_n(K)^F = \mathrm{PGL}_n(q)$ where $n = 2^r \times k$, both $r$ and $k$ are positive integers, and $q$ is a power of $\ell$. Let $\Lambda$ be a nontoral elementary abelian 2-subgroup of $\mathrm{PGL}_n(K)^F = \mathrm{PGL}_n(q)$. Let $E \leq \mathrm{PGL}_n(\mathbb{C})$ correspond to $\Lambda$ under the bijection of Proposition 3.1. Now
$$E = \bar{\Gamma}_r \times \bar{A}$$

where $\bar{A}$ is either trivial or a representative of a conjugacy class of the toral elementary abelian 2-subgroups of $\text{PGL}_k(\mathbb{C})$. By Proposition 3.2, $\Lambda = 2^{2r} \times \Lambda_w$ where $\Lambda_w$ is either trivial or a toral elementary abelian 2-subgroup of $\text{PGL}_k(K)$. We denote $N_{\text{PGL}_n(K)}(\Lambda)$ by $N$ and $C_{\text{PGL}_n(K)}(\Lambda)$ by $C$. We study the $F$-action on $N/C^\circ$.

If $\bar{A}$ is a maximal toral elementary abelian 2-subgroup of $\text{PGL}_k(\mathbb{C})$, then $\bar{\Gamma}_r \times \bar{A}$ is a maximal nontoral elementary abelian 2-subgroup of $\text{PGL}_n(\mathbb{C})$ and this is already considered. So we study those nonmaximal toral $\bar{A}$. Let $Q = \bar{\Gamma}_r \times J$ be a maximal nontoral elementary abelian 2-subgroup of $\text{PGL}_n(\mathbb{C})$ where $J$ is a maximal toral elementary abelian 2-subgroup of $\text{PGL}_k(\mathbb{C})$ of rank $m$ where $m = k-1$. Let $E = \bar{\Gamma}_r \times L$ be a nontoral elementary abelian 2-subgroup of $\text{PGL}_n(\mathbb{C})$ where $L < J$. Now let $\Lambda = 2^{2r} \times \Lambda_w$ correspond to $E$. Here

$$\begin{aligned} C_{\text{PGL}_n(\mathbb{C})}(E) &\cong \bar{\Gamma}_r \times C_{\text{PGL}_k(\mathbb{C})}(L) \\ C_{\text{PGL}_n(\mathbb{C})}(Q) &\cong \bar{\Gamma}_r \times T_m \\ W_{\text{PGL}_n(\mathbb{C})}(E) &= 2^{\text{rk}L \times 2r} : (\text{Sp}_{2r}(2) \times W_{\text{PGL}_k(\mathbb{C})}(L)) \\ W_{\text{PGL}_n(\mathbb{C})}(Q) &= 2^{m \times 2r} : (\text{Sp}_{2r}(2) \times S_{m+1}). \end{aligned}$$

We denote by $U_w$ the group $2^{\text{rk}L \times 2r}$ in $W_{\text{PGL}_n(\mathbb{C})}(E)$.

If $r \geq 1$, then $\bar{\Gamma}_r.\text{Sp}_{2r}(2)$ is a subgroup of $\text{PGL}_{2^r}(q)$ for $q$ a power of $\ell$ and thus, $\bar{\Gamma}_r.\text{Sp}_{2r}(2)$ is centralized by $F$. By Lemma 4.5, we conclude that $F$ also centralizes $(C_{\text{PGL}_k(K)}(\Lambda_w)/C_{\text{PGL}_k(K)}(\Lambda_w)^\circ).W_{\text{PGL}_k(K)}(\Lambda_w)$. So our focus is on the $U_w$ part in $W_{\text{PGL}_n(\mathbb{C})}(E)$.

Recall the groups $R$ and $N_1$ from Section 4.2.

**Theorem 4.6.** *Let $G = \text{PGL}_n(\mathbb{C})$ where $n = 2^r \times k$ and both $r$ and $k$ are positive integers. Let $Q = \bar{\Gamma}_r \times J$ be a maximal nontoral elementary abelian 2-subgroup of $G$ where $J$ is a maximal toral elementary abelian 2-subgroup of $\text{PGL}_k(\mathbb{C})$ of rank $m$ where $m = k - 1$. Let $E = \bar{\Gamma}_r \times L$ be a nontoral elementary abelian 2-subgroup of $G$ where $L < J$. There exists $U_1 \leq N_1$ such that $u \in N_G(E) \setminus C_G(E)$ for each nontrivial $u \in U_1$ and $U_1 C_G(E)/C_G(E)$ is the subgroup $U_w$ of $W_G(E)$.*

*Proof.* For $1 \leq i \leq m$ and $1 \leq j \leq r$, the generators of $N_1$ are $x_{ij}$ and $y_{ij}$ where

$$x_{ij} = \left[ \begin{pmatrix} \Delta(I_{2^r})_{i-1} & & \\ & \Delta(\begin{smallmatrix} & -I_{2^{j-1}} \\ I_{2^{j-1}} & \end{smallmatrix})_{2^{r-j}} & \\ & & \Delta(I_{2^r})_{m-(i-1)} \end{pmatrix} \right],$$

$$y_{ij} = \left[ \begin{pmatrix} \Delta(I_{2^r})_{i-1} & & \\ & \Delta(\begin{smallmatrix} -I_{2^{j-1}} & \\ & I_{2^{j-1}} \end{smallmatrix})_{2^{r-j}} & \\ & & \Delta(I_{2^r})_{m-(i-1)} \end{pmatrix} \right]$$

and $x_{ij}, y_{ij} \in N_G(Q) \setminus C_G(Q)$.

The group $J$ has generators $g_i = [M_i]$ where $i = 1,...,m$ and $M_i$ is a $k \times k$ diagonal matrix with $-1$ in the $i$-th position and $1$ in the remaining positions. Hence $Q$ has generators $\bar{A}_s$, $\bar{B}_s$ and $g_i \otimes I_{2^r}$ for $0 \leq s \leq r-1$ and $i = 1,...,m$. Computation shows

$$x_{ij}\bar{A}_s x_{ij}^{-1} = \begin{cases} \bar{A}_s & s+1 \neq j \\ \bar{A}_s(g_i \otimes I_{2^r}) & s+1 = j \end{cases}$$

$$x_{ij}\bar{B}_s x_{ij}^{-1} = \begin{cases} \bar{B}_s & s+1 \neq j \\ \bar{B}_s(g_i \otimes I_{2^r}) & s+1 = j \end{cases}$$

$$y_{ij}\bar{A}_s y_{ij}^{-1} = \bar{A}_s \text{ for all } s$$

$$y_{ij}\bar{B}_s y_{ij}^{-1} = \begin{cases} \bar{B}_s & s+1 \neq j \\ \bar{B}_s(g_i \otimes I_{2^r}) & s+1 = j \end{cases}$$

$x_{ij}$ and $y_{ij}$ fix all the $g_i \otimes I_{2^r}$ for $i = 1,...,m$.

Let the generators of $L$ be $h_1,...,h_t$ where $t < m$. So $E$ has generators $\bar{A}_s$, $\bar{B}_s$ and $h_i \otimes I_{2^r}$ for $0 \leq s \leq r-1$ and $i = 1,...,t$. Since $L < J$, each of the generators of $L$ is a product of some generators of $J$. For instance, let $h_1 = g_{z_1}g_{z_2}...g_{z_r}$. Correspondingly, we form products $v_{1j} = x_{z_1 j}\, x_{z_2 j}... x_{z_r j}$ and $w_{1j} = y_{z_1 j}\, y_{z_2 j}... y_{z_r j}$ for $1 \leq j \leq r$. In general, corresponding to the $t$ generators of $L$, we form $2rt$ products $v_{ij}$ and $w_{ij}$ for $1 \leq i \leq t$ and $1 \leq j \leq r$.

The $v_{ij}$ and $w_{ij}$ generate an elementary abelian 2-subgroup $U_1$ of rank $2rt$. Also,

$$v_{ij}\bar{A}_s v_{ij}^{-1} = \begin{cases} \bar{A}_s & s+1 \neq j \\ \bar{A}_s(h_i \otimes I_{2^r}) & s+1 = j \end{cases}$$

$$v_{ij}\bar{B}_s v_{ij}^{-1} = \begin{cases} \bar{B}_s & s+1 \neq j \\ \bar{B}_s(h_i \otimes I_{2^r}) & s+1 = j \end{cases}$$

$$w_{ij}\bar{A}_s w_{ij}^{-1} = \bar{A}_s \text{ for all } s$$

$$w_{ij}\bar{B}_s w_{ij}^{-1} = \begin{cases} \bar{B}_s & s+1 \neq j \\ \bar{B}_s(h_i \otimes I_{2^r}) & s+1 = j \end{cases}$$

$v_{ij}$ and $w_{ij}$ fix all the $h_i \otimes I_{2^r}$ for $\beta = 1,...,t$.

Clearly $U_1 \leq N_1$ and $u \in N_G(E) \setminus C_G(E)$ for every nontrivial $u \in U_1$. By the action of $U_1$ on $E$ and the structure of $W_G(E)$, we conclude that $U_1 C_G(E)/C_G(E)$ is the subgroup $U_w$ of $W_G(E)$. $\square$

Consequently, as in the maximal nontoral case, $N/C^\circ$ is centralized by $F$ when $\Lambda$ is a nonmaximal nontoral elementary abelian 2-subgroup of $\mathrm{PGL}_n(K)$ and $q$ is a power of a prime $\ell$ where $\ell \equiv 1 \bmod 4$.

4.3.2. *The (local) structure of toral 2-subgroups with disconnected centralizers.*

Corresponding to $E = \bar{\Gamma}_r \times \bar{A}$ is $\Lambda = 2^{2r} \times \Lambda_w$. We first study the structure of $\Lambda_w$ with a disconnected centralizer in $\mathrm{PGL}_k(K)$. Then we examine its local structure in Theorems 4.10 and 4.14.

**Theorem 4.7** ([18, Theorem 14.2]). *Let $G$ be a connected reductive algebraic group, let $s \in G$ be semisimple, and let $T \leq G$ be a maximal torus with corresponding root system $\Phi$. Let $s \in T$ and $\Psi = \{\alpha \in \Phi \,|\, \alpha(s) = 1\}$. Let $\dot{w}$ be a choice of preimage in $N_G(T)$ of $w \in W$ and let $s_\alpha$ be a reflection along $\alpha$.*

$$\begin{aligned} C_G(s)^\circ &= \langle T, U_\alpha \,|\, \alpha \in \Psi \rangle \\ C_G(s) &= \langle T, U_\alpha, \dot{w} \,|\, \alpha \in \Psi, w \in W \text{ with } s^w = s \rangle. \end{aligned}$$

*Moreover, $C_G(s)^\circ$ is reductive with root system $\Psi$ and Weyl group $W_1 = \langle s_\alpha \,|\, \alpha \in \Psi \rangle$.*

**Remark 4.8.** *Let $s$ and $t$ be semisimple elements of $T$, $\Psi_1 := \{\alpha \in \Phi \,|\, \alpha(s) = 1\}$ and $\Psi_2 := \{\beta \in \Phi \,|\, \beta(t) = 1\}$. By Theorem 4.7:*

(i) *If $C_G(s) = C_G(s)^\circ = \langle T, U_\alpha \,|\, \alpha \in \Psi_1 \rangle$ and $C_G(t) = C_G(t)^\circ = \langle T, U_\beta \,|\, \beta \in \Psi_2 \rangle$, then $C_G(\langle s,t \rangle) = \langle T, U_\gamma \,|\, \gamma \in \Psi_1 \bigcap \Psi_2 \rangle$, which is connected. This result can be generalized to finitely many semisimple elements in $T$ with connected centralizers.*

(ii) *If*
$$C_G(s) = \langle T, U_\alpha, \dot{w}_1 \,|\, \alpha \in \Psi_1, w_1 \in W \text{ with } s^{w_1} = s \rangle$$
*and*
$$C_G(t) = \langle T, U_\beta, \dot{w}_2 \,|\, \beta \in \Psi_2, w_2 \in W \text{ with } t^{w_2} = t \rangle,$$
*then $C_G(\langle s,t \rangle) = \langle T, U_\gamma, X, Y \,|\, \gamma \in \Psi_1 \bigcap \Psi_2 \rangle$ where $X := \langle \dot{w}_1 \rangle \bigcap C_G(t)$ and $Y := \langle \dot{w}_2 \rangle \bigcap C_G(s)$. This result can be generalized to finitely many semisimple elements in $T$.*

(iii) *If $E_1, E_2 \leq T$ are toral elementary abelian subgroups of $G$, then*
$$C_G(E_1)^\circ \bigcap C_G(E_2)^\circ = C_G(\langle E_1, E_2 \rangle)^\circ.$$

Let $G = \mathrm{PGL}_n(K)$ where $n = 2^s \times t$, $\gcd(2,t) = 1$ and $s \geq 1$. Without loss of generality, let the maximal torus $T$ be the image in $G$ of the group of invertible diagonal matrices in $\mathrm{GL}_n(K)$. We choose the representatives in $T$ of the conjugacy classes of the toral elementary abelian 2-subgroups of $G$, so the elements of the representatives are the images in $G$ of the diagonal matrices with diagonal entries 1 or $-1$.

As in Section 4.1, $\bar{\Gamma}_r = \langle \bar{A}_0, \bar{B}_0, \bar{A}_1, \bar{B}_1, ..., \bar{A}_{r-1}, \bar{B}_{r-1} \rangle$. Let $D_r = \langle \bar{A}_0, \bar{A}_1, ..., \bar{A}_{r-1} \rangle$ and $B_r = \langle \bar{B}_0, \bar{B}_1, ..., \bar{B}_{r-1} \rangle$. Note that $D_r$ and $B_r$ can be viewed as subgroups of $G$.

**Theorem 4.9.** *Let $n = 2^s \times t$ with $\gcd(2,t) = 1$ and $s, t \geq 1$. The conjugacy classes of the toral elementary abelian 2-subgroups of $\mathrm{PGL}_n(K)$ with disconnected centralizers have representatives $D_r \times H$ for $1 \leq r \leq s$ where $n = 2^r \times k$, and $H$ is either trivial or a representative of a conjugacy class of the toral elementary abelian 2-subgroups of $\mathrm{PGL}_k(K)$ with $C_{\mathrm{PGL}_k(K)}(H)$ connected.*

*Proof.* Let $G = \mathrm{PGL}_n(K)$. We use induction. When $n = 2$, it is true. Assume that the conclusion holds for $n = 2, 4, \ldots, 2u$ where $u \geq 1$.

Let $n = 2(u+1)$ and let $G_1 = \mathrm{PGL}_{u+1}(K)$. By [14, Table 4.3.1], the only conjugacy class of involutions in $G$ with disconnected centralizers has representative

$$e_{\frac{n}{2}} = \left[ \begin{pmatrix} -I_{\frac{n}{2}} & \\ & I_{\frac{n}{2}} \end{pmatrix} \right]$$

and

$$C_G(e_{\frac{n}{2}}) = C_G(e_{\frac{n}{2}})^\circ \sqcup f_{\frac{n}{2}} C_G(e_{\frac{n}{2}})^\circ,$$

viewed as a disjoint union of two cosets where

$$f_{\frac{n}{2}} = \left[ \begin{pmatrix} & I_{\frac{n}{2}} \\ I_{\frac{n}{2}} & \end{pmatrix} \right].$$

For $v \geq 1$, let $D = \langle g_1, g_2, \ldots, g_v \rangle \leqslant T$ be a toral elementary abelian 2-subgroup of $G$ with $C_G(D) \neq C_G(D)^\circ$.

By Remark 4.8 (i), if no element of $D$ is conjugate to $e_{\frac{n}{2}}$, then $C_G(D)$, as the intersection of the centralizers of all elements of $D$, is connected. So there exists at least one element of $D$ conjugate to $e_{\frac{n}{2}}$. For every $e' \in D$ conjugate to $e_{\frac{n}{2}}$, there exists $f' \in G$ such that $C_G(e') = \langle C_G(e')^\circ, f' \rangle$. By Remark 4.8 (ii), if none of the $f'$ centralizes $D$, then $C_G(D)$ is connected. Therefore, there exists $e' \in D$ conjugate to $e_{\frac{n}{2}}$ with

$$C_G(e') = \langle C_G(e')^\circ, f' \rangle$$

and $f'$ centralizes every element of $D$. Moreover, by conjugating $D$, we may identify $e'$ with $e_{\frac{n}{2}}$, so $f' = f_{\frac{n}{2}}$. We may also assume $g_v = e_{\frac{n}{2}}$. Since $f_{\frac{n}{2}}$ centralizes each of the remaining generators, the form which $g_1, g_2, \ldots, g_{v-1}$ take is $[\Delta(S)_2]$ where

$$S = \begin{pmatrix} * & & \\ & * & \\ & & \ddots \\ & & & * \end{pmatrix}_{(u+1) \times (u+1)} \quad \text{and the diagonal entries are 1 or } -1.$$

Hence $\langle g_1, g_2, \ldots, g_{v-1} \rangle$ is either trivial or isomorphic to a toral elementary abelian 2-subgroup $H$ of $G_1$. We may further conjugate $D$ such that $g_v$ becomes $\bar{A}_0$ and $f_{\frac{n}{2}}$ becomes $\bar{B}_0$. So

$$D = \langle \bar{A}_0 \rangle \times H = D_1 \times H.$$

If $u+1$ is odd, then $C_{G_1}(H)$ is connected for each toral elementary abelian 2-subgroup $H$ (or the trivial subgroup) of $G_1$ and the conjugacy classes of the toral elementary abelian 2-subgroups of $G$ with disconnected centralizers have representatives $D_1 \times H$ where $H$ is either trivial or a representative of a conjugacy class of the toral elementary abelian 2-subgroups of $G_1$.

If $u + 1 = 2^{s-1} \times t$ is even and $H$ is either trivial or a representative of a conjugacy class of the toral elementary abelian 2-subgroups of $G_1$ with $C_{G_1}(H)$ connected, then $D_1 \times H$ is a representative of a conjugacy class of the toral elementary abelian 2-subgroups of $G$ with disconnected centralizers.

Now assume $u + 1 = 2^{s-1} \times t$ is even and $H$ is a representative of a conjugacy class of the toral elementary abelian 2-subgroups of $G_1$ with $C_{G_1}(H) \neq C_{G_1}(H)^\circ$. By induction, for each factorization $u + 1 = 2^w \times k$ where $1 \leq w \leq s - 1$,
$$H = \left\langle \bar{A}'_0, \bar{A}'_1, ..., \bar{A}'_{w-1} \right\rangle \times J$$
where
$$\bar{A}'_j = \left[ \Delta\left( \begin{pmatrix} -I_{2^j} & \\ & I_{2^j} \end{pmatrix} \right)_{2^{s-1-j-1}t} \right]$$
for $0 \leq j \leq w - 1$, and $J$ is either trivial or a representative of a conjugacy class of the toral elementary abelian 2-subgroups of $\mathrm{PGL}_k(K)$ with $C_{\mathrm{PGL}_k(K)}(J)$ connected. Now $n = 2^s \times t = 2^{w+1} \times k$ and
$$\begin{aligned} D &= D_1 \times (H \otimes I_2) \\ &= D_1 \times ((\langle \bar{A}'_0, \bar{A}'_1, ..., \bar{A}'_{w-1}\rangle \times (J \otimes I_{2^w})) \otimes I_2) \\ &= \langle \bar{A}_0, \bar{A}_1, ..., \bar{A}_w \rangle \times (J \otimes I_{2^{w+1}}) \\ &= D_{w+1} \times (J \otimes I_{2^{w+1}}) \end{aligned}$$
where $J$ is either trivial or a representative of a conjugacy class of the toral elementary abelian 2-subgroups of $\mathrm{PGL}_k(K)$ with $C_{\mathrm{PGL}_k(K)}(J)$ connected. $\square$

Next, we obtain the structure of the centralizer quotient.

**Theorem 4.10.** *Let $n = 2^s \times t$ with $\gcd(2, t) = 1$ and $s, t \geq 1$. For $n = 2^r \times k$ where $1 \leq r \leq s$, if $D = D_r \times H$ in $G = \mathrm{PGL}_n(K)$ where $H$ is either trivial or a representative of a conjugacy class of the toral elementary abelian 2-subgroups of $\mathrm{PGL}_k(K)$ with $C_{\mathrm{PGL}_k(K)}(H)$ connected, then*
$$C_G(D) = C_G(D)^\circ.B_r.$$

*Proof.* First, $C_G\left(\bar{A}_i\right) = \left\langle C_G(\bar{A}_i)^\circ, \bar{B}_i \right\rangle$ for $i = 0, ..., r - 1$. Each nontrivial element of $\bar{\Gamma}_r$ is conjugate to $e_{\frac{n}{2}}$, and each $\bar{B}_i$ centralizes $\bar{A}_0, \bar{A}_1, ..., \bar{A}_{r-1}$. By Remark 4.8 (ii),
$$C_G(D_r) = \left\langle C_G(D_r)^\circ, \bar{B}_i \,|\, i = 0, ..., r - 1 \right\rangle.$$

Let $h_j$ for $j = 1, ..., w$ be the generators of $H$. If $h_j \otimes I_{2^r}$ is not conjugate to $\bar{A}_0$ for $j = 1, ..., w$, then $C_G(h_j \otimes I_{2^r})$ is connected and
$$C_G(D) = \langle C_G(D)^\circ, \bar{B}_i \,|\, i = 0, ..., r-1 \rangle,$$
by Remark 4.8 (ii), because all of the $\bar{B}_i$ centralize $D$.

Now suppose $h_j \otimes I_{2^r}$ is conjugate to $\bar{A}_0$ for some $1 \leq j \leq w$. We can conjugate $H$ in $\mathrm{PGL}_k(K)$ such that the conjugate of $h_j$, which we denote again by $h_j$, is now
$$\left[ \begin{pmatrix} -I_{k/2} & \\ & I_{k/2} \end{pmatrix} \right]$$
and $C_G(h_j \otimes I_{2^r}) = \langle C_G(h_j \otimes I_{2^r})^\circ, y \rangle$ where $y = y_j \otimes I_{2^r}$ and
$$y_j = \left[ \begin{pmatrix} & I_{k/2} \\ I_{k/2} & \end{pmatrix} \right].$$
We conclude that $y$ cannot centralize $H \otimes I_{2^r}$: otherwise, $y_j$ would centralize $H$ and this would imply that $y_j \in C_{\mathrm{PGL}_k(K)}(H) \setminus C_{\mathrm{PGL}_k(K)}(H)^\circ$. But $C_{\mathrm{PGL}_k(K)}(H)$ is connected. So we have a contradiction. Since each $\bar{B}_i$ centralizes $D$,
$$C_G(D) = \langle C_G(D)^\circ, \bar{B}_i \,|\, i = 0, ..., r-1 \rangle.$$
Thus, the theorem is proved. $\square$

**Theorem 4.11.** *Let $n = 2^s \times t$ with $\gcd(2, t) = 1$ and $s, t \geq 1$. If $D$ is a toral elementary abelian 2-subgroup of $G = \mathrm{PGL}_n(K)$ with a disconnected centralizer, then*
$$N_G(D)/C_G(D) \geq \mathrm{GL}_r(2)$$
*where $1 \leq r \leq s$ and $n = 2^r \times 2^{s-r} t$. Hence there are two $F$-classes of $N_G(D)/C_G(D)^\circ$ in $C_G(D)/C_G(D)^\circ$, one consisting of the identity and the other the remaining elements.*

*Proof.* By Theorem 4.9, $D = D_r \times H \leq \Lambda \cong \bar{\Gamma}_r \times H$. Now $C_G(\Lambda) \cong \bar{\Gamma}_r \times C_{\mathrm{PGL}_{2^{s-r}t}(K)}(H)$,
$$W_G(\Lambda) = 2^{\mathrm{rk}H \times 2r} : (\mathrm{Sp}_{2r}(2) \times W_{\mathrm{PGL}_{2^{s-r}t}(K)}(H))$$
and
$$D.(\mathrm{GL}_r(2) \times \mathrm{Id}) \leq \Lambda.(\mathrm{Sp}_{2r}(2) \times \mathrm{Id}) \leq N_G(\Lambda).$$
Let $\Omega$ be an index set such that for $g_i \in N_G(\Lambda)$ where $i \in \Omega$,
$$\langle g_i C_G(\Lambda) \,|\, i \in \Omega \rangle / C_G(\Lambda) \cong \mathrm{GL}_r(2) \times \mathrm{Id}.$$
So
$$M := \langle g_i C_G(\Lambda) \,|\, i \in \Omega \rangle \leq N_G(D).$$
The centralizer
$$C_G(D) \leq \langle g_i C_G(D) \,|\, i \in \Omega \rangle.$$

So
$$Q := \langle g_i C_G(D) \,|\, i \in \Omega \rangle \leq N_G(D).$$
Let $f = f_2 \circ f_1$ where $f_1 : M \to Q$ is the natural embedding and $f_2 : Q \to Q/C_G(D)$ is the natural surjective map. Also, $f$ maps $g_i c \in M$ where $c \in C_G(\Lambda)$ to $g_i C_G(D)$ and $f$ is surjective. Now $\ker f = M \cap C_G(D)$. We prove that $\ker f = C_G(\Lambda)$. Clearly $M \cap C_G(D) \geq C_G(\Lambda)$. Consider the other direction. Suppose there is some $gc \in M$ where $g \notin C_G(\Lambda)$ and $c \in C_G(\Lambda) \leq C_G(D)$ such that $gc \in C_G(D)$. Since $M/C_G(\Lambda) \cong \mathrm{GL}_r(2) \times \mathrm{Id}$ acts faithfully on $D$, we deduce that $g \in C_G(\Lambda)$, a contradiction. Hence
$$\mathrm{GL}_r(2) \times \mathrm{Id} \cong M/C_G(\Lambda) \cong Q/C_G(D) \leq N_G(D)/C_G(D).$$
By Theorem 4.10, $C_G(D) = C_G(D)^\circ . B_r$ where $B_r = \langle \bar{B}_0, \bar{B}_1, ..., \bar{B}_{r-1} \rangle$. Now
$$\mathrm{GL}_r(2) \cong \left\{ \begin{pmatrix} Y & 0 \\ 0 & Y^{-T} \end{pmatrix} \leq \mathrm{Sp}_{2r}(2) \,\Big|\, Y \in \mathrm{GL}_r(2) \right\}$$
where $-T$ is transpose inverse.

Let $U = 2^{1+2r}$ be extraspecial and $\bar{U} = 2^{2r}$ the quotient of $U$ by its center $\langle z \rangle$. For $x \in U$, denote by $\bar{x}$ its image in $\bar{U}$. By [21, Section 2], if we set $\beta(\bar{x}, \bar{y}) = \bar{c}$ where $[x,y] = z^c$ for $x, y \in U$, then $\beta$ is a nondegenerate alternating bilinear form on $\bar{U}$. Now we identify $\{\bar{A}_0, \bar{B}_0, \bar{A}_1, \bar{B}_1, ..., \bar{A}_{r-1}, \bar{B}_{r-1}\}$ with the corresponding vector space and equip it with $\beta$. Thus, we can view $\bar{\Gamma}_r$ as a symplectic vector space over $\mathrm{GF}(2)$ and $\{\bar{A}_0, \bar{B}_0, \bar{A}_1, \bar{B}_1, ..., \bar{A}_{r-1}, \bar{B}_{r-1}\}$ as a basis of this space. By the action of $\mathrm{GL}_r(2)$ on the basis, we conclude that there are two $F$-classes in $C_{\mathrm{PGL}_n(K)}(D)/C_{\mathrm{PGL}_n(K)}(D)^\circ$, one consisting of the identity and the other the remaining elements. $\square$

**Theorem 4.12.** *Let $G = \mathrm{PGL}_n(K)$ where $n \geq 2$ is even. Let $D$ be a representative of a conjugacy class of toral elementary abelian 2-subgroups of $G$. If $C_G(D) = C_G(D)^\circ$, then $D$ has an involution not conjugate to*
$$e_{\frac{n}{2}} = \left[ \begin{pmatrix} -I_{\frac{n}{2}} & \\ & I_{\frac{n}{2}} \end{pmatrix} \right].$$

*Proof.* We use induction on the rank of $D$ to show this. When $D$ has rank 1, it is true. When $D$ has rank 2, we suppose that the three involutions in $D$ are conjugate to $e_{\frac{n}{2}}$ and one is $e_{\frac{n}{2}}$. Now
$$C_G(e_{\frac{n}{2}}) = C_G(e_{\frac{n}{2}})^\circ \sqcup f_{\frac{n}{2}} C_G(e_{\frac{n}{2}})^\circ$$
where
$$C_G(e_{\frac{n}{2}})^\circ = \mathrm{SL}_{\frac{n}{2}}(K) \circ_{(n/2)^*} T_1 \circ_{(n/2)^*} \mathrm{SL}_{\frac{n}{2}}(K) \text{ and } f_{\frac{n}{2}} = \left[ \begin{pmatrix} & I_{\frac{n}{2}} \\ I_{\frac{n}{2}} & \end{pmatrix} \right].$$

(Recall that $(n/2)^*$ denotes $\gcd(n/2, \ell-1)$ where $\ell$ is the characteristic of $K$.) Since $D$ is toral, $D \leq C_G(e_{\frac{n}{2}})^\circ$. Let $e'$ be another involution of $D$ which has the shape

$$\left[\begin{pmatrix} e_u & \\ & e_d \end{pmatrix}\right]$$

where $e_u$ and $e_d$ are diagonal matrices from $\mathrm{GL}_{\frac{n}{2}}(K)$. Let $x$ be the number of $-1$ on the diagonal of $e_u$ and let $y$ be the number of $-1$ on the diagonal of $e_d$. Since $e' \sim e_{\frac{n}{2}}$ and $e'e_{\frac{n}{2}} \sim e_{\frac{n}{2}}$, we deduce that

(1) $$\begin{cases} x + y = \frac{n}{2} \\ (\frac{n}{2} - x) + y = \frac{n}{2} \end{cases}$$

which yields

the number of $-1$ of $e_u =$ the number of $-1$ of $e_d = n/4$.

Consequently, we can conjugate $D$ such that

$$e_u = e_d = \begin{pmatrix} -I_{\frac{n}{4}} & \\ & I_{\frac{n}{4}} \end{pmatrix}.$$

Hence $f_{\frac{n}{2}}$ centralizes $e'$ and therefore $D$. So $C_{\mathrm{PGL}_n(K)}(D)$ is disconnected, a contradiction.

Now we assume the statement is true when $D$ has rank $k$.

Assume $D$ has rank $k+1$. If it has an elementary abelian 2-subgroup $E$ of rank $k$ and $C_{\mathrm{PGL}_n(K)}(E)$ is connected, then by assumption $E$ has an involution not conjugate to $e_{\frac{n}{2}}$. Now consider the case where none of the elementary abelian subgroups of rank $k$ of $D$ has a connected centralizer in $G$. Let $n = 2^s \times t, s \geq 1, \gcd(2,t) = 1$ and assume $D = \langle g_1, g_2, ..., g_{k+1} \rangle$. Let $E_1 = \langle g_1, g_2, ..., g_k \rangle$ and $E_2 = \langle g_2, ..., g_{k+1} \rangle$. Since both $C_{\mathrm{PGL}_n(K)}(E_1)$ and $C_{\mathrm{PGL}_n(K)}(E_2)$ are disconnected, by Theorem 4.9, $E_1$ and $E_2$ take the shape $D_r \times H$ for some $1 \leq r \leq s, n = 2^r \times k$ and $H$ is either trivial or a representative of a conjugacy class of the toral elementary abelian 2-subgroups of $\mathrm{PGL}_k(K)$ with $C_{\mathrm{PGL}_k(K)}(H)$ connected. In particular, by Theorem 4.10, $\bar{B}_0$ centralizes both $E_1$ and $E_2$, but $\bar{B}_0$ is not in $C_{\mathrm{PGL}_n(K)}(E_1)^\circ$ or $C_{\mathrm{PGL}_n(K)}(E_2)^\circ$. Hence $\bar{B}_0$ centralizes $D$ but

$$\bar{B}_0 \notin C_{\mathrm{PGL}_n(K)}(D)^\circ = C_{\mathrm{PGL}_n(K)}(E_1)^\circ \cap C_{\mathrm{PGL}_n(K)}(E_2)^\circ.$$

So $C_{\mathrm{PGL}_n(K)}(D)$ is disconnected, a contradiction. □

**Theorem 4.13.** *Let $G = \mathrm{PGL}_n(K)$ where $n \geq 2$ is even. Let $D$ be a representative of a conjugacy class of toral elementary abelian 2-subgroups of $G$. If $C_G(D) = C_G(D)^\circ$,* *then there exists a choice of generators of $D$ such that none is conjugate to*

$$e_{\frac{n}{2}} = \left[\begin{pmatrix} -I_{\frac{n}{2}} & \\ & I_{\frac{n}{2}} \end{pmatrix}\right].$$

*Proof.* We use induction on the rank of $D$. If $D$ has rank 1, then the claim holds: clearly, there is no toral elementary abelian 2-subgroup of $\mathrm{PGL}_2(\mathbb{C})$ with connected centralizer. If $D$ has rank 2, then by Theorem 4.12 there is an involution $e'$ not conjugate to $e_{\frac{n}{2}}$. Let

$$e' = \left[\begin{pmatrix} -I_k & \\ & I_{n-k} \end{pmatrix}\right]$$

where $1 \leq k < n/2$ and let

$$D = \langle e', e_1 \rangle.$$

If $e_1$ is not conjugate to $e_{\frac{n}{2}}$, then the claim is true. Next, suppose $e_1$ is conjugate to $e_{\frac{n}{2}}$. If $e'e_1$ is not conjugate to $e_{\frac{n}{2}}$, then the claim is true once we replace $e_1$ by $e'e_1$. Now consider the case where $e'e_1$ is conjugate to $e_{\frac{n}{2}}$. Since

$$C_G(e') = \mathrm{SL}_k(K) \circ_{k^*} T_1 \circ_{(n-k)^*} \mathrm{SL}_{n-k}(K),$$

each involution of $D$ takes the shape

$$\left[\begin{pmatrix} M_1 & \\ & M_2 \end{pmatrix}\right]$$

where $M_1$ is a $k \times k$ diagonal matrix, $M_2$ is a $(n-k) \times (n-k)$ diagonal matrix and both have 1 or $-1$ as diagonal entries. Let $x$ denote the number of $-1$ on the diagonal of $M_1$ in $e_1$ and $y$ the number of $-1$ on the diagonal of $M_2$ in $e_1$. Since $e'e_1$ is conjugate to $e_{\frac{n}{2}}$, we deduce that

(2) $$\begin{cases} x + y = \frac{n}{2} \\ (k-x) + y = \frac{n}{2} \end{cases}$$

which yields a unique solution

(3) $$\begin{cases} x = \frac{k}{2} \\ y = \frac{n-k}{2}. \end{cases}$$

Thus,

$$e_1 = \left[\begin{pmatrix} \begin{pmatrix} -I_{\frac{k}{2}} & \\ & I_{\frac{k}{2}} \end{pmatrix} & \\ & \begin{pmatrix} -I_{\frac{n-k}{2}} & \\ & I_{\frac{n-k}{2}} \end{pmatrix} \end{pmatrix}\right]$$

and $D = \langle e', e_1 \rangle$ has a disconnected centralizer in $G$, a contradiction.

Assume the claim is true when $D$ has rank $k \geq 2$.

Now consider the case where $D$ has rank $k + 1$. First, we suppose there exists a subgroup $E$ of $D$ of rank $k$ with a connected centralizer in $G$. By assumption, there is a choice of generators $g_1, g_2, \ldots, g_k$ of $E$ such that none is conjugate to $e_{\frac{n}{2}}$.

Let $D = \langle g_1, g_2, \ldots, g_k, g_{k+1} \rangle$. If $g_{k+1}$ is not conjugate to $e_{\frac{n}{2}}$, then the claim is true. If $g_{k+1}$ is conjugate to $e_{\frac{n}{2}}$ but $g_{k+1}e$ is not conjugate to $e_{\frac{n}{2}}$ for some $e \in E$, then the claim is true once we replace $g_{k+1}$ by $g_{k+1}e$. Suppose $g_{k+1}$ is conjugate to $e_{\frac{n}{2}}$ and $g_{k+1}e$ is conjugate to $e_{\frac{n}{2}}$ for all $e \in E$. Without loss of generality, let

$$g_1 = \left[ \begin{pmatrix} -I_{k_1} & \\ & I_{n-k_1} \end{pmatrix} \right]$$

where $1 \leq k_1 < n/2$. Since

$$C_G(g_1) = \mathrm{SL}_{k_1}(K) \circ_{(k_1)^*} T_1 \circ_{(n-k_1)^*} \mathrm{SL}_{n-k_1}(K),$$

each involution of $D$ takes the shape

$$\left[ \begin{pmatrix} M_1 & \\ & M_2 \end{pmatrix} \right]$$

where $M_1$ is a $k_1 \times k_1$ diagonal matrix, $M_2$ is a $(n - k_1) \times (n - k_1)$ diagonal matrix and both have 1 or $-1$ as diagonal entries. Since $g_1 g_{k+1}$ is conjugate to $e_{\frac{n}{2}}$ and $g_{k+1}$ is conjugate to $e_{\frac{n}{2}}$, by the argument for rank 2,

$$g_{k+1} = \left[ \begin{pmatrix} \begin{pmatrix} -I_{\frac{k_1}{2}} & \\ & I_{\frac{k_1}{2}} \end{pmatrix} & \\ & \begin{pmatrix} -I_{\frac{n-k_1}{2}} & \\ & I_{\frac{n-k_1}{2}} \end{pmatrix} \end{pmatrix} \right].$$

Now let

$$g_2 = \left[ \begin{pmatrix} \begin{pmatrix} -I_{k_{21}} & \\ & I_{k_1 - k_{21}} \end{pmatrix} & \\ & \begin{pmatrix} -I_{k_{22}} & \\ & I_{n-k_1-k_{22}} \end{pmatrix} \end{pmatrix} \right].$$

First, $k_{21} + k_{22} < n/2$ since $g_2$ is not conjugate to $e_{\frac{n}{2}}$. There are five possibilities for $k_{21}$ and $k_{22}$:

(1) $k_{21} < \frac{k_1}{2}$, $k_{22} < \frac{n-k_1}{2}$
(2) $k_{21} = \frac{k_1}{2}$, $k_{22} < \frac{n-k_1}{2}$
(3) $k_{21} < \frac{k_1}{2}$, $k_{22} = \frac{n-k_1}{2}$
(4) $k_{21} > \frac{k_1}{2}$, $k_{22} < \frac{n-k_1}{2}$
(5) $k_{21} < \frac{k_1}{2}$, $k_{22} > \frac{n-k_1}{2}$.

Hence we can assume without loss of generality that $k_{21} \leq k_1/2, k_{22} \leq (n-k_1)/2$ but equality does not hold simultaneously: if possibility 4 or 5 holds, then we can multiply $g_2$ by

$$g_1 = \left[\begin{pmatrix} -I_{k_1} & \\ & I_{n-k_1} \end{pmatrix}\right] = \left[\begin{pmatrix} I_{k_1} & \\ & -I_{n-k_1} \end{pmatrix}\right]$$

and $g_1 g_2$ is not conjugate to $e_{\frac{n}{2}}$. Since $g_2 g_{k+1}$ is conjugate to $e_{\frac{n}{2}}$,

$$\frac{k_1}{2} - k_{21} + \left(\frac{n-k_1}{2}\right) - k_{22} = \frac{n}{2}.$$

This yields $-k_{21} = k_{22}$, which is impossible unless $k_{21} = k_{22} = 0$.

Finally, we assume there is no subgroup $E$ of $D$ of rank $k$ with a connected centralizer in $G$. By an argument similar to the last paragraph in the proof of Theorem 4.12, we reach a contradiction. □

Now we give the precise structure of the normalizer quotient.

**Theorem 4.14.** *Let $G = \mathrm{PGL}_n(K)$ where $n = 2^s \times t$ with $\gcd(2, t) = 1$ and $s, t \geq 1$. For $n = 2^r \times 2^{s-r} t$ where $1 \leq r \leq s$, if $D = D_r \times \bar{U}$ where $\bar{U} = \bar{U}_1 \otimes I_{2^r}$ and $\bar{U}_1$ is either trivial or a representative of a conjugacy class of toral elementary abelian 2-subgroups of $\mathrm{PGL}_{2^{s-r} t}(K)$ with $C_{\mathrm{PGL}_{2^{s-r} t}(K)}(\bar{U}_1)$ connected, then*

$$N_G(D)/C_G(D) \cong 2^{\mathrm{rk} \bar{U} \times r} : (\mathrm{GL}_r(2) \times W_{\mathrm{PGL}_{2^{s-r} t}(K)}(\bar{U}_1)).$$

*Proof.* Let $Z$ be the center of $\mathrm{GL}_n(K)$. Let $D \leq E = \bar{\Gamma}_r \times \bar{U}$. So

$$W_G(E) = 2^{\mathrm{rk} \bar{U} \times 2r} : (\mathrm{Sp}_{2r}(2) \times W_{\mathrm{PGL}_{2^{s-r} t}(K)}(\bar{U}_1)).$$

Therefore,

$$N_G(D)/C_G(D) \leq \begin{pmatrix} \mathrm{GL}_r(2) & 0 \\ 2^{\mathrm{rk} \bar{U} \times r} & \mathrm{GL}_{\mathrm{rk} \bar{U}}(2) \end{pmatrix} = 2^{\mathrm{rk} \bar{U} \times r} : (\mathrm{GL}_r(2) \times \mathrm{GL}_{\mathrm{rk} \bar{U}}(2)).$$

We denote the subgroup $2^{\mathrm{rk} \bar{U} \times r}$ by $R$.

For $g \in N_G(D)$, let $\bar{g} = \begin{pmatrix} g_1 & 0 \\ g_2 & g_3 \end{pmatrix}$ where $g_1 \in \mathrm{GL}_r(2)$, $g_2 \in R$ and $g_3 \in \mathrm{GL}_{\mathrm{rk} \bar{U}}(2)$. Since there exist elements

$$\bar{h} = \begin{pmatrix} g_1^{-1} & 0 \\ 0 & I_{\mathrm{rk} \bar{U}} \end{pmatrix} \text{ and } \bar{y} = \begin{pmatrix} I_r & 0 \\ -g_2 & I_{\mathrm{rk} \bar{U}} \end{pmatrix}$$

of $N_G(D)/C_G(D)$ such that $\bar{y} \bar{h} \bar{g} = \begin{pmatrix} I_r & 0 \\ 0 & g_3 \end{pmatrix}$, we can assume $g \in N_{C_G(D_r)}(\bar{U})$.

We are interested in $N_{C_G(D_r)}(\bar{U})/C_G(D) = N_{C_G(D_r)}(\bar{U})/C_{C_G(D_r)}(\bar{U})$. Observe that
$$C_G(D_r) = ((\mathrm{GL}_{2^{s-r}t}(K) \times \ldots \times \mathrm{GL}_{2^{s-r}t}(K))/Z).B_r$$
where the direct product has $2^r$ copies of $\mathrm{GL}_{2^{s-r}t}(K)$. Let $U_1 \leq \mathrm{GL}_{2^{s-r}t}(K)$ with $U_1 Z/Z = \bar{U}_1$. Now
$$N_{C_G(D_r)}(\bar{U}) \leq (N_{\mathrm{GL}_{2^{s-r}t}(K)}(U_1 Z))^{2^r}/Z,$$
so
$$\frac{N_{C_G(D_r)}(\bar{U})}{C_{C_G(D_r)}(\bar{U})} \leq \frac{(N_{\mathrm{GL}_{2^{s-r}t}(K)}(U_1 Z))^{2^r}/Z}{C_{C_G(D_r)}(\bar{U})}.$$
If $g_j$ is from the $j$-th copy of $N_{\mathrm{GL}_{2^{s-r}t}(K)}(U_1 Z)$ for $2 \leq j \leq 2^r$, then there exists $g_j^{-1}$ in the first copy of $N_{\mathrm{GL}_{2^{s-r}t}(K)}(U_1 Z)$ such that in $(N_{\mathrm{GL}_{2^{s-r}t}(K)}(U_1 Z))^{2^r}/Z$
$$(1, \ldots, g_j, \ldots, 1) Z = (g_j, 1, \ldots, 1)(g_j^{-1}, \ldots, g_j, \ldots, 1) Z.$$
Moreover,
$$(g_j^{-1}, \ldots, g_j, \ldots, 1) Z \in C_{C_G(D_r)}(\bar{U}).$$
Therefore,
$$\frac{(N_{\mathrm{GL}_{2^{s-r}t}(K)}(U_1 Z))^{2^r}/Z}{C_{C_G(D_r)}(\bar{U})} \cong \frac{N_{\mathrm{GL}_{2^{s-r}t}(K)}(U_1 Z)/Z}{C_{C_G(D_r)}(\bar{U}) \cap N_{\mathrm{GL}_{2^{s-r}t}(K)}(U_1 Z)/Z}$$
$$\cong \frac{N_{\mathrm{PGL}_{2^{s-r}t}(K)}(\bar{U}_1)}{C_{C_G(D_r)}(\bar{U}) \cap N_{\mathrm{PGL}_{2^{s-r}t}(K)}(\bar{U}_1)}$$
$$= \frac{N_{\mathrm{PGL}_{2^{s-r}t}(K)}(\bar{U}_1)}{C_{\mathrm{PGL}_{2^{s-r}t}(K)}(\bar{U}_1)}$$
and
$$\frac{N_{C_G(D_r)}(\bar{U})}{C_{C_G(D_r)}(\bar{U})} \leq \frac{N_{\mathrm{PGL}_{2^{s-r}t}(K)}(\bar{U}_1)}{C_{\mathrm{PGL}_{2^{s-r}t}(K)}(\bar{U}_1)}.$$
On the other hand, if $y \in N_{\mathrm{GL}_{2^{s-r}t}(K)}(U_1 Z)$, then $(y \otimes I_{2^r})Z/Z \in N_{C_G(D_r)}(\bar{U})$.

Hence the other inclusion holds. $\square$

4.3.3. *Descending to the finite groups.*

Let $\Lambda$ be a nontoral elementary abelian 2-subgroup of $\mathrm{PGL}_n(K)^F = \mathrm{PGL}_n(q)$ where $q \equiv 1 \mod 4$ and let $E \leq \mathrm{PGL}_n(\mathbb{C})$ correspond to $\Lambda$ under the bijection of Proposition 3.1. Now $E = \bar{\Gamma}_r \times \bar{A}$ where $\bar{A}$ is either trivial or a representative of a conjugacy class of the toral elementary abelian 2-subgroups of $\mathrm{PGL}_k(\mathbb{C})$. We denote $N_{\mathrm{PGL}_n(K)}(\Lambda)$ by $N$ and $C_{\mathrm{PGL}_n(K)}(\Lambda)$ by $C$.

If $\bar{A}$ is trivial, then $E = \bar{\Gamma}_r$, and $C_{\mathrm{PGL}_n(\mathbb{C})}(E) \cong \bar{\Gamma}_r \times \mathrm{PGL}_k(\mathbb{C})$ and $W_{\mathrm{PGL}_n(\mathbb{C})}(E) =$

$\mathrm{Sp}_{2r}(2)$. For $r \geq 1$, the group $\bar{\Gamma}_r.\mathrm{Sp}_{2r}(2)$ is a subgroup of $\mathrm{PGL}_{2^r}(q)$ for $q$ a power of a prime $\ell$ where $\ell \equiv 1 \bmod 4$ and is centralized by $F$. Hence, when $\Lambda \cong 2^{2r}$, there are two $F$-classes of $N/C^\circ$ in $C/C^\circ$, one consisting of the identity and the other all the remaining elements. Correspondingly, there are two $\mathrm{PGL}_n(q)$-conjugacy classes of elementary abelian 2-subgroups with representatives $E_1$ and $E_2$. Their local structure is the following:

$$\begin{aligned} C_{\mathrm{PGL}_n(q)}(E_1) &\cong \bar{\Gamma}_r \times \mathrm{PGL}_k(q) \\ N_{\mathrm{PGL}_n(q)}(E_1) &\cong (\bar{\Gamma}_r \times \mathrm{PGL}_k(q)).\mathrm{Sp}_{2r}(2) \\ C_{\mathrm{PGL}_n(q)}(E_2) &\cong \bar{\Gamma}_r \times \mathrm{PGL}_k(q) \\ N_{\mathrm{PGL}_n(q)}(E_2) &\cong (\bar{\Gamma}_r \times \mathrm{PGL}_k(q)).(2^{2r-1} : \mathrm{Sp}_{2(r-1)}(2)). \end{aligned}$$

Now we focus on nontrivial $\bar{A}$ and discuss separately the cases where $\bar{A}$ has a connected and disconnected centralizer.

**Case 1**. $\bar{A}$ is nontrivial nonmaximal toral in $\mathrm{PGL}_k(\mathbb{C})$ with $C_{\mathrm{PGL}_k(\mathbb{C})}(\bar{A})$ connected.

Now
$$C_{\mathrm{PGL}_n(\mathbb{C})}(E) \cong \bar{\Gamma}_r \times C_{\mathrm{PGL}_k(\mathbb{C})}(\bar{A}).$$
Since $C_{\mathrm{PGL}_k(\mathbb{C})}(\bar{A})$ is connected, $C_{\mathrm{PGL}_n(\mathbb{C})}(E)/C_{\mathrm{PGL}_n(\mathbb{C})}(E)^\circ \cong \bar{\Gamma}_r$. By Proposition 3.2, $\Lambda = 2^{2r} \times \Lambda_w$ where $\Lambda_w$ is either trivial or a toral elementary abelian 2-subgroup of $\mathrm{PGL}_k(K)$ with a connected centralizer. When $\Lambda$ is a nonmaximal nontoral elementary abelian 2-subgroup of $\mathrm{PGL}_n(K)$, $N/C^\circ$ is centralized by $F$.

By construction of the generators of $U_1$ from Theorem 4.6, $C_{\mathrm{PGL}_n(\mathbb{C})}(E)/C_{\mathrm{PGL}_n(\mathbb{C})}(E)^\circ$ is centralized by $U_1$. Hence, when $\Lambda = 2^{2r} \times \Lambda_w$ where $\Lambda_w$ has a connected centralizer in $\mathrm{PGL}_k(K)$, there are two $F$-classes of $N/C^\circ$ in $C/C^\circ$, one consisting of the identity and the other all of the remaining elements. Correspondingly, there are two $\mathrm{PGL}_n(q)$-conjugacy classes of elementary abelian 2-subgroups with representatives $E_1$ and $E_2$. Their local structure is the following:

$$\begin{aligned} C_{\mathrm{PGL}_n(q)}(E_1) &\cong \bar{\Gamma}_r \times \left(C_{\mathrm{PGL}_k(K)}(\Lambda_w)\right)^F \\ C_{\mathrm{PGL}_n(q)}(E_2) &\cong \bar{\Gamma}_r \times \left(C_{\mathrm{PGL}_k(K)}(\Lambda_w)\right)^F \\ N_{\mathrm{PGL}_n(q)}(E_1) &\cong C_{\mathrm{PGL}_n(q)}(E_1).W_1 \\ N_{\mathrm{PGL}_n(q)}(E_2) &\cong C_{\mathrm{PGL}_n(q)}(E_2).W_2 \end{aligned}$$

where
$$\begin{aligned} W_1 &= 2^{\mathrm{rk}\bar{A}\times 2r} : (\mathrm{Sp}_{2r}(2) \times W_{\mathrm{PGL}_k(\mathbb{C})}(\bar{A})) \\ W_2 &= 2^{\mathrm{rk}\bar{A}\times 2r} : ((2^{2r-1} : \mathrm{Sp}_{2(r-1)}(2)) \times W_{\mathrm{PGL}_k(\mathbb{C})}(\bar{A})). \end{aligned}$$

**Case 2**. $\bar{A}$ is nonmaximal toral in $\mathrm{PGL}_k(\mathbb{C})$ with $C_{\mathrm{PGL}_k(\mathbb{C})}(\bar{A})$ disconnected.

Let $G = \mathrm{PGL}_n(K)$ where $n = 2^s \times t$, $\gcd(2,t) = 1$ and $s \geq 1$. For each factorization

$n = 2^r \times k$ where $1 \leq r \leq s$, let $\Lambda \cong \bar{\Gamma}_r \times \Lambda_w$ be a representative of a conjugacy class of the nonmaximal nontoral elementary abelian 2-subgroups where $\Lambda_w$ is a representative of a conjugacy class of the nonmaximal toral elementary abelian 2-subgroups of $\mathrm{PGL}_k(K)$ with $C_{\mathrm{PGL}_k(K)}(\Lambda_w)$ disconnected.

By Theorem 4.9, for each factorization $k = 2^{r_1} \times k_1$ where both $r_1$ and $k_1$ are positive integers, $\Lambda_w = D_{r_1} \times H$. Here $D_{r_1} = \langle \bar{A}_0^k, ..., \bar{A}_{r_1-1}^k \rangle$ where the superscript $k$ means elements in $\mathrm{PGL}_k(K)$ and $H$ is either trivial or a representative of a conjugacy class of the toral elementary abelian 2-subgroups of $\mathrm{PGL}_{k_1}(K)$ with $C_{\mathrm{PGL}_{k_1}(K)}(H)$ connected. Let $B_{r_1} = \langle \bar{B}_0^k, ..., \bar{B}_{r_1-1}^k \rangle$. Now

$$\begin{aligned} C_G(\Lambda) &\cong \bar{\Gamma}_r \times (C_{\mathrm{PGL}_k(K)}(\Lambda_w)^\circ . B_{r_1}) \\ N_G(\Lambda)/C_G(\Lambda) &\cong 2^{\mathrm{rk}\Lambda_w \times 2r} : (\mathrm{Sp}_{2r}(2) \times W_{\mathrm{PGL}_k(K)}(\Lambda_w)) \end{aligned}$$

where

$$W_{\mathrm{PGL}_k(K)}(\Lambda_w) \cong 2^{\mathrm{rk}H \times r_1} : (\mathrm{GL}_{r_1}(2) \times W_{\mathrm{PGL}_{k_1}(K)}(H)).$$

We denote by $R_1$ the group $2^{\mathrm{rk}\Lambda_w \times 2r}$ in $N_G(\Lambda)/C_G(\Lambda)$.

The element $\bar{A}_0^k$ of $D_{r_1}$ is a product of some generators of the maximal toral elementary abelian 2-subgroup of $\mathrm{PGL}_k(K)$. Corresponding to this product, we construct $y \in R_1$ in the same way as in the proof of Theorem 4.6 to get:

$$y = \left[ \Delta\left( \begin{pmatrix} \Delta((^{-1}\ 1))_{2^{r-1}} & \\ & I_{2^r} \end{pmatrix} \right)_{\frac{k}{2}} \right] C_G(\Lambda).$$

Now

$$\bar{B}_0^k = \left[ \Delta\left( \begin{pmatrix} & 1 \\ 1 & \end{pmatrix} \right)_{\frac{k}{2}} \right].$$

Thus,

$$y(\bar{B}_0^k \otimes I_{2^r}) y^{-1} = (\bar{B}_0^k \otimes I_{2^r}) \bar{A}_0.$$

By Theorem 4.11 and the action of $R_1$ on $\bar{\Gamma}_r$, we conclude that when $\Lambda = 2^{2r} \times \Lambda_w$ where $\Lambda_w$ has a disconnected centralizer in $\mathrm{PGL}_k(K)$, there are three $F$-classes of $N/C^\circ$ in $C/C^\circ$. One consists of the identity, another consists of the nontrivial elements of $2^{2r}$ and the third consists of all the remaining elements with representative $b$. Correspondingly, there are three $\mathrm{PGL}_n(q)$-conjugacy classes of elementary abelian

2-subgroups with representatives $E_1$, $E_2$ and $E_3$. Their local structure is the following:

$$\begin{aligned} C_{\mathrm{PGL}_n(q)}(E_1) &\cong \bar{\Gamma}_r \times ((C_{\mathrm{PGL}_k(K)}(\Lambda_w)^\circ)^F . B_{r_1}) \\ C_{\mathrm{PGL}_n(q)}(E_2) &\cong \bar{\Gamma}_r \times ((C_{\mathrm{PGL}_k(K)}(\Lambda_w)^\circ)^F . B_{r_1}) \\ C_{\mathrm{PGL}_n(q)}(E_3) &\cong \bar{\Gamma}_r \times ((C_{\mathrm{PGL}_k(K)}(\Lambda_w)^\circ)^{bF} . B_{r_1}) \\ N_{\mathrm{PGL}_n(q)}(E_1) &\cong C_{\mathrm{PGL}_n(q)}(E_1).(2^{\mathrm{rk}\Lambda_w \times 2r} : (\mathrm{Sp}_{2r}(2) \times W_{\mathrm{PGL}_k(K)}(\Lambda_w))) \\ N_{\mathrm{PGL}_n(q)}(E_2) &\cong C_{\mathrm{PGL}_n(q)}(E_2).(2^{\mathrm{rk}\Lambda_w \times 2r} : ((2^{2r-1} : \mathrm{Sp}_{2(r-1)}(2)) \times W_{\mathrm{PGL}_k(K)}(\Lambda_w))) \\ N_{\mathrm{PGL}_n(q)}(E_3) &\cong C_{\mathrm{PGL}_n(q)}(E_3).(2^{(\mathrm{rk}\Lambda_w-1)\times 2r} : (\mathrm{Sp}_{2r}(2) \times W_3)) \end{aligned}$$

where

$$W_3 \cong 2^{\mathrm{rk}H \times r_1} : (\begin{pmatrix} 1 & * \\ 0 & \mathrm{GL}_{r_1-1}(2) \end{pmatrix} \times W_{\mathrm{PGL}_{k_1}(K)}(H)).$$

This finishes the discussion of the case where $\ell \equiv 1 \bmod 4$. Next, we look at the case where $\ell \equiv 3 \bmod 4$.

### 4.4. Nontoral 2-subgroups of $\mathrm{PGL}_n(q)$ where $q \equiv 3 \bmod 4$

Let $K$ be an algebraically closed field of characteristic $\ell$ where $\ell \equiv 3 \bmod 4$. Let $F$ be a Steinberg endomorphism such that $\mathrm{PGL}_n(K)^F = \mathrm{PGL}_n(q)$ where $q$ is a power of $\ell$. In this section, we discuss the elementary abelian 2-subgroups of $\mathrm{PGL}_n(q)$. Here the $F$-action need not be trivial; we consider it in Theorem 4.17. Then we descend to the finite groups and discuss separately the cases where $\bar{A}$ has a connected and disconnected centralizer.

First, we suppose $q$ is an odd power of $\ell$. Below we record a result used to investigate the nontrivial $F$-action.

**Proposition 4.15** ([4, Proposition 4.5]). *Let $G$ be a connected reductive algebraic group over an algebraically closed field of characteristic $\ell$. Let $F$ be a Steinberg endomorphism of $G$. Let $L$ be a finite subgroup of $G^F$. The map $N_G(L) \to N_G(L)$ induced by $F$ is given by $n \mapsto n\theta(n)$, where $\theta$ is a 1-cocycle $N_G(L) \to C_G(L)$; that is, $\theta$ satisfies*

$$\theta(n_1 n_2) = (\theta(n_1)^{n_2})\theta(n_2).$$

*In particular, $F$ induces the identity map on $N_G(L)/C_G(L)$.*

Let $\Lambda$ be a nontoral elementary abelian 2-subgroup of $\mathrm{PGL}_n(K)^F = \mathrm{PGL}_n(q)$. Let $E \leq \mathrm{PGL}_n(\mathbb{C})$ correspond to $\Lambda$ under the bijection of Proposition 3.1. Now

$$E = \bar{\Gamma}_r \times \bar{A}$$

where $\bar{A}$ is either trivial or a representative of a conjugacy class of the toral elementary abelian 2-subgroups of $\mathrm{PGL}_k(\mathbb{C})$. By Proposition 3.2, $\Lambda = 2^{2r} \times \Lambda_w$ where $\Lambda_w$ is either trivial or a toral elementary abelian 2-subgroup of $\mathrm{PGL}_k(K)$. We denote $N_{\mathrm{PGL}_n(K)}(\Lambda)$ by $N$ and $C_{\mathrm{PGL}_n(K)}(\Lambda)$ by $C$. Now
$$C_{\mathrm{PGL}_n(\mathbb{C})}(E) \cong \bar{\Gamma}_r \times C_{\mathrm{PGL}_k(\mathbb{C})}(\bar{A})$$
and
$$W_{\mathrm{PGL}_n(\mathbb{C})}(E) = 2^{\mathrm{rk}\bar{A} \times 2r} : (\mathrm{Sp}_{2r}(2) \times W_{\mathrm{PGL}_k(\mathbb{C})}(\bar{A})).$$
We denote the subgroup $2^{\mathrm{rk}\bar{A} \times 2r}$ by $R$.

By Proposition 3.4, to find the classes of subgroups of $\mathrm{PGL}_n(q)$ which are $\mathrm{PGL}_n(K)$-conjugate to $\Lambda$, we must consider the $F$-classes of $N/C^\circ$ in $C/C^\circ$.

By Lemma 4.5, $(C_{\mathrm{PGL}_k(K)}(\Lambda_w)/C_{\mathrm{PGL}_k(K)}(\Lambda_w)^\circ).W_{\mathrm{PGL}_k(K)}(\Lambda_w)$ is centralized by $F$. By the construction in Section 4.2 of the generators of $R$, we deduce that $F$ also centralizes $R$. So we focus on the $F$-action on $\bar{\Gamma}_r.\mathrm{Sp}_{2r}(2)$.

**Lemma 4.16** ([18, Corollary 21.8]). *Let $G$ be a connected reductive algebraic group with a Steinberg endomorphism $F : G \to G$. In the semidirect product $G : \langle F \rangle$, the coset $G.F$ of $F$ consists of a single conjugacy class, so $G.F = F^G$ where $F^G$ is the orbit of $F$ under the action of $G$. In particular, $(G^g)^F$ and $G^F$ are $G$-conjugate for every $g \in G$.*

Let $E = \bar{\Gamma}_r$ be an elementary abelian 2-subgroup of $\mathrm{PGL}_{2^r}(\mathbb{C})$ of order $2^{2r}$. Clearly
$$C_{\mathrm{PGL}_{2^r}(\mathbb{C})}(E) \cong \bar{\Gamma}_r \text{ and } N_{\mathrm{PGL}_{2^r}(\mathbb{C})}(E) \cong \bar{\Gamma}_r.\mathrm{Sp}_{2r}(2).$$
Now we investigate the nontrivial $F$-action.

**Theorem 4.17.** *Let $H = \mathrm{PGL}_{2^r}(K)$ where $r \geq 1$ and $K$ is an algebraically closed field of characteristic $\ell$ where $\ell \equiv 3 \bmod 4$. Let $E$ be a nontoral elementary abelian subgroup of $H$ of order $2^{2r}$. Let $F$ be a Steinberg endomorphism such that $H^F = \mathrm{PGL}_{2^r}(q)$ where $q$ is an odd power of $\ell$. There exist two $F$-classes of $N_H(E)/C_H(E)^\circ$ in $C_H(E)/C_H(E)^\circ$. Correspondingly, there are two $H^F$-conjugacy classes of elementary abelian 2-subgroups of order $2^{2r}$ with representatives $E_1$ and $E_2$. Their local structure is the following:*
$$C_{H^F}(E_1) \cong 2^{2r} \text{ and } N_{H^F}(E_1) \cong 2^{2r}.\mathrm{SO}_{2r}^+(2)$$
$$C_{H^F}(E_2) \cong 2^{2r} \text{ and } N_{H^F}(E_2) \cong 2^{2r}.\mathrm{SO}_{2r}^-(2).$$

*Proof.* Let $U = \mathrm{GL}_{2^r}(q)$ where $q$ is an odd power of $\ell$ and $\ell \equiv 3 \bmod 4$. Denote the center of $U$ by $C$.

By [2, 1E], there are extraspecial subgroups $S_1 := 2^{1+2r}_+$ and $S_2 := 2^{1+2r}_-$ of $U$ with $N_U(S_1) = CS_1.\mathrm{SO}_{2r}^+(2)$ and $N_U(S_2) = CS_2.\mathrm{SO}_{2r}^-(2)$. We deduce that
$$2^{2r}.\mathrm{SO}_{2r}^+(2) \leq H^F$$

and
$$2^{2r}.\mathrm{SO}_{2r}^{-}(2) \le H^F.$$
Moreover,
$$2^{2r}.\mathrm{SO}_{2r}^{+}(2) \le 2^{2r}.\mathrm{Sp}_{2r}(2)$$
and
$$2^{2r}.\mathrm{SO}_{2r}^{-}(2) \le 2^{2r}.\mathrm{Sp}_{2r}(2)$$
since $\mathrm{SO}_{2r}^{+}(2)$ and $\mathrm{SO}_{2r}^{-}(2)$ are maximal subgroups of $\mathrm{Sp}_{2r}(2)$ by [10, Section 8.2] and [17, Proposition 4.8.6]. Thus, there are at least two $F$-classes of $N_H(E)/C_H(E)^\circ$ in $C_H(E)/C_H(E)^\circ$.

On the other hand, suppose $2^{2r}.\mathrm{SO}_{2r}^{+}(2)$ is centralized by $wF$ for some $w \in 2^{2r}$. Lemma 4.16 allows us to replace $F$ by $wF$. Thus, we can assume that $2^{2r}.\mathrm{SO}_{2r}^{+}(2)$ is centralized by $F$. Observe that $2^{2r}.\mathrm{SO}_{2r}^{+}(2) \le 2^{2r}.\mathrm{Sp}_{2r}(2)$. We identify $2^{2r}$ with the corresponding vector space and equip it with the quadratic form provided by [17, Section 2]. Hence, by [17, Lemma 2.10.5], the action of $2^{2r}.\mathrm{SO}_{2r}^{+}(2)$ on $2^{2r}$ yields orbits $O_1$, $O_2$ and $O_3$ where $O_1$ contains the zero vector, $O_2$ contains the isotropic vectors, and $O_3$ contains the nonisotropic vectors.

Observe that $F$ acts nontrivially on $2^{2r}.\mathrm{Sp}_{2r}(2)$: otherwise, $q$ is a power of a prime $\ell$ where $\ell \equiv 1 \bmod 4$. But $F$ acts trivially on $N_H(E)/C_H(E)$ by Proposition 4.15. Thus, there exists $n \in 2^{2r}.\mathrm{Sp}_{2r}(2)$ such that $F(n)n^{-1} \ne \mathrm{Id}$. Hence $F$ fuses $O_1$ with one of $O_2$ or $O_3$ and there are at most two $F$-classes of $N_H(E)/C_H(E)^\circ$ in $C_H(E)/C_H(E)^\circ$.

Consequently, there are exactly two $F$-classes of $N_H(E)/C_H(E)^\circ$ in $C_H(E)/C_H(E)^\circ$. Correspondingly, there are two $H^F$-conjugacy classes of elementary abelian subgroups of order $2^{2r}$ with representatives $E_1$ and $E_2$. Their local structure is the following:
$$C_{H^F}(E_1) \cong 2^{2r} \text{ and } N_{H^F}(E_1) \cong 2^{2r}.\mathrm{SO}_{2r}^{+}(2)$$
$$C_{H^F}(E_2) \cong 2^{2r} \text{ and } N_{H^F}(E_2) \cong 2^{2r}.\mathrm{SO}_{2r}^{-}(2).$$
Thus, the theorem is proved. $\square$

### Descending to the finite groups

When $\bar{A}$ is trivial, $E = \bar{\Gamma}_r$, $C_{\mathrm{PGL}_n(\mathbb{C})}(E) \cong \bar{\Gamma}_r \times \mathrm{PGL}_k(\mathbb{C})$ and $W_{\mathrm{PGL}_n(\mathbb{C})}(E) = \mathrm{Sp}_{2r}(2)$. Hence, by Theorem 4.17, there are two $F$-classes of $N/C^\circ$ in $C/C^\circ$ when $\Lambda \cong 2^{2r}$, and, correspondingly, there are two $\mathrm{PGL}_n(q)$-conjugacy classes of elementary abelian 2-subgroups with representatives $E_1$ and $E_2$. Their local structure is the following:
$$C_{\mathrm{PGL}_n(q)}(E_1) \cong \bar{\Gamma}_r \times \mathrm{PGL}_k(q) \text{ and } N_{\mathrm{PGL}_n(q)}(E_1) \cong (\bar{\Gamma}_r \times \mathrm{PGL}_k(q)).\mathrm{SO}_{2r}^{+}(2)$$
$$C_{\mathrm{PGL}_n(q)}(E_2) \cong \bar{\Gamma}_r \times \mathrm{PGL}_k(q) \text{ and } N_{\mathrm{PGL}_n(q)}(E_2) \cong (\bar{\Gamma}_r \times \mathrm{PGL}_k(q)).\mathrm{SO}_{2r}^{-}(2).$$

We focus on the nontrivial $\bar{A}$ and consider two cases.

**Case 1**. $\bar{A}$ is nontrivial toral in $\mathrm{PGL}_k(\mathbb{C})$ with $C_{\mathrm{PGL}_k(\mathbb{C})}(\bar{A})$ connected.

As a representative of a conjugacy class of the nontoral elementary abelian 2-subgroups of $\mathrm{PGL}_n(\mathbb{C})$,
$$E = \bar{\Gamma}_r \times \bar{A}$$
where $\bar{A}$ is a representative of a conjugacy class of the nontrivial toral elementary abelian 2-subgroups of $\mathrm{PGL}_k(\mathbb{C})$ and $C_{\mathrm{PGL}_k(\mathbb{C})}(\bar{A})$ is connected.

Since $C_{\mathrm{PGL}_k(\mathbb{C})}(\bar{A})$ is connected, $C_{\mathrm{PGL}_n(\mathbb{C})}(E)/C_{\mathrm{PGL}_n(\mathbb{C})}(E)^\circ \cong \bar{\Gamma}_r$. By construction of the generators of $U_1$ from Theorem 4.6, $C_{\mathrm{PGL}_n(\mathbb{C})}(E)/C_{\mathrm{PGL}_n(\mathbb{C})}(E)^\circ$ is centralized by $U_1$. By Theorem 4.17, there are two $F$-classes of $N/C^\circ$ in $C/C^\circ$ when $\Lambda = 2^{2r} \times \Lambda_w$ where $\Lambda_w$ has a connected centralizer in $\mathrm{PGL}_k(K)$. Correspondingly, there are two $\mathrm{PGL}_n(q)$-conjugacy classes of elementary abelian 2-subgroups with representatives $E_1$ and $E_2$. Their local structure is the following:

$$\begin{aligned} C_{\mathrm{PGL}_n(q)}(E_1) &\cong \bar{\Gamma}_r \times \left(C_{\mathrm{PGL}_k(K)}(\Lambda_w)\right)^F \\ C_{\mathrm{PGL}_n(q)}(E_2) &\cong \bar{\Gamma}_r \times \left(C_{\mathrm{PGL}_k(K)}(\Lambda_w)\right)^F \\ N_{\mathrm{PGL}_n(q)}(E_1) &\cong C_{\mathrm{PGL}_n(q)}(E_1).(2^{\mathrm{rk}\bar{A}\times 2r} : (\mathrm{SO}_{2r}^+(2) \times W_{\mathrm{PGL}_k(\mathbb{C})}(\bar{A}))) \\ N_{\mathrm{PGL}_n(q)}(E_2) &\cong C_{\mathrm{PGL}_n(q)}(E_2).(2^{\mathrm{rk}\bar{A}\times 2r} : (\mathrm{SO}_{2r}^-(2) \times W_{\mathrm{PGL}_k(\mathbb{C})}(\bar{A}))). \end{aligned}$$

**Case 2**. $\bar{A}$ is toral in $\mathrm{PGL}_k(\mathbb{C})$ with $C_{\mathrm{PGL}_k(\mathbb{C})}(\bar{A})$ disconnected.

Let $\Lambda = 2^{2r} \times \Lambda_w$ be a representative of a conjugacy class of the nontoral elementary abelian 2-subgroups of $G = \mathrm{PGL}_n(K)$ where $\Lambda_w$ is a representative of a conjugacy class of the toral elementary abelian 2-subgroups of $\mathrm{PGL}_k(K)$ with $C_{\mathrm{PGL}_k(K)}(\Lambda_w)$ disconnected.

For each factorization $k = 2^{r_1} \times k_1$ where both $r_1$ and $k_1$ are positive integers, $\Lambda_w = D_{r_1} \times H$. Here $D_{r_1} = \langle \bar{A}_0^k, ..., \bar{A}_{r-1}^k \rangle$ where the superscript $k$ means elements in $\mathrm{PGL}_k(K)$, and $H$ is either trivial or a representative of a conjugacy class of the toral elementary abelian 2-subgroups of $\mathrm{PGL}_{k_1}(K)$ with $C_{\mathrm{PGL}_{k_1}(K)}(H)$ connected. Let $B_{r_1} = \langle \bar{B}_0^k, ..., \bar{B}_{r-1}^k \rangle$. Now

$$\begin{aligned} C_G(\Lambda) &\cong \bar{\Gamma}_r \times (C_{\mathrm{PGL}_k(K)}(\Lambda_w)^\circ . B_{r_1}) \\ N_G(\Lambda)/C_G(\Lambda) &\cong R_1 : (\mathrm{Sp}_{2r}(2) \times W_{\mathrm{PGL}_k(K)}(\Lambda_w)). \end{aligned}$$

where $R_1 = 2^{\mathrm{rk}\Lambda_w \times 2r}$. So

$$C_G(\Lambda)/C_G(\Lambda)^\circ \cong \bar{\Gamma}_r \times B_{r_1}.$$

Recall from Case 2 of Section 4.3.3 that when $q$ is a power of a prime $\ell$ where $\ell \equiv 1 \bmod 4$, there are two $F$-classes $O_1, O_2$ in $\bar{\Gamma}_r$ resulting from the $F$-action of $\mathrm{Sp}_{2r}(2)$ with $O_1$ containing the identity only and $O_2$ the remaining elements. The $F$-action of $\mathrm{GL}_{r_1}(2)$ on $B_{r_1}$ gives two $F$-classes: $O_3$ contains the identity, and $O_4$ consists of all other elements. Let $Or_1 = \{xy \,|\, x \in O_1, y \in O_4\}$ and $Or_2 = \{xy \,|\, x \in O_2, y \in O_4\}$. When $q$ is a power of a prime $\ell$ where $\ell \equiv 1 \mod 4$, we deduce that $Or_1$ and $Or_2$ are fused into one class by the action of $R_1$.

Assume that $q$ is an odd power of a prime $\ell$ where $\ell \equiv 3 \mod 4$. By Theorem 4.17, there are two $F$-classes $O'_1$, $O'_2$ in $\bar{\Gamma}_r$ resulting from the $F$-action of $\mathrm{Sp}_{2r}(2)$ with $O'_1$ containing the identity and some other elements. The two $F$-classes $O_3$ and $O_4$ in $B_{r_1}$ resulting from the $F$-action of $\mathrm{GL}_{r_1}(2)$ are the same as the case where $\ell \equiv 1 \mod 4$. Let $Or'_1 = \{xy \,|\, x \in O'_1, y \in O_4\}$ and $Or'_2 = \{xy \,|\, x \in O'_2, y \in O_4\}$. Note that $Or_1 \subset Or'_1$ and $Or'_2 \subset Or_2$. Let $u \in Or_1$ and $v \in Or'_2$. Since $F$ acts trivially on $R_1$, we use the results from the previous paragraph to deduce that there exists $r \in R_1$ such that $rur^{-1} = v$. Thus, $Or'_1$ and $Or'_2$ are fused. Therefore, there are three $F$-classes of $N/C^\circ$ in $C/C^\circ$ when $\Lambda = 2^{2r} \times \Lambda_w$ where $\Lambda_w$ has a disconnected centralizer in $\mathrm{PGL}_k(K)$. The class $Or'_1 \cup Or'_2$ has representative $b'$. Correspondingly, there are three $\mathrm{PGL}_n(q)$-conjugacy classes of elementary abelian 2-subgroups with representatives $E_1$, $E_2$ and $E_3$. Their local structure is the following:

$$
\begin{aligned}
C_{\mathrm{PGL}_n(q)}(E_1) &\cong \bar{\Gamma}_r \times ((C_{\mathrm{PGL}_k(K)}(\Lambda_w)^\circ)^F.B_{r_1}) \\
C_{\mathrm{PGL}_n(q)}(E_2) &\cong \bar{\Gamma}_r \times ((C_{\mathrm{PGL}_k(K)}(\Lambda_w)^\circ)^F.B_{r_1}) \\
C_{\mathrm{PGL}_n(q)}(E_3) &\cong \bar{\Gamma}_r \times ((C_{\mathrm{PGL}_k(K)}(\Lambda_w)^\circ)^{b'F}.B_{r_1}) \\
N_{\mathrm{PGL}_n(q)}(E_1) &\cong C_{\mathrm{PGL}_n(q)}(E_1).(2^{\mathrm{rk}\Lambda_w \times 2r} : (\mathrm{SO}^+_{2r}(2) \times W_{\mathrm{PGL}_k(K)}(\Lambda_w))) \\
N_{\mathrm{PGL}_n(q)}(E_2) &\cong C_{\mathrm{PGL}_n(q)}(E_2).(2^{\mathrm{rk}\Lambda_w \times 2r} : (\mathrm{SO}^-_{2r}(2) \times W_{\mathrm{PGL}_k(K)}(\Lambda_w))) \\
N_{\mathrm{PGL}_n(q)}(E_3) &\cong C_{\mathrm{PGL}_n(q)}(E_3).(2^{(\mathrm{rk}\Lambda_w - 1) \times 2r} : (\mathrm{Sp}_{2r}(2) \times W_3))
\end{aligned}
$$

where

$$
W_3 \cong 2^{\mathrm{rk}H \times r_1} : (\begin{pmatrix} 1 & * \\ 0 & \mathrm{GL}_{r_1-1}(2) \end{pmatrix} \times W_{\mathrm{PGL}_{k_1}(K)}(H)).
$$

Now we consider the case where $q$ is an even power of $\ell$ where $\ell \equiv 3 \mod 4$. We adopt some notation from [1, Section 1] and record two results. Let $\mathrm{GF}(q)$ have odd characteristic and let $2^{a+1}$ be the exact power of 2 dividing $q^2 - 1$. So $a \geq 2$. The prime 2 is *linear* or *unitary* according as $2^a$ divides $q - 1$ or $q + 1$. If $p$ is a prime, then a $p$-group is of *symplectic-type* if every characteristic abelian subgroup is cyclic.

**Proposition 4.18** ([1, Proposition 1A]). *Let $E$ be the extraspecial group of order $2^{2r+1}$ with type $\eta$, where $r \geq 1$ and $\eta = +$ or $-$. There exists a unique faithful and irreducible representation of $E$ in $\mathrm{GL}_{2^r}(q)$.*

**Proposition 4.19** ([1, Proposition 1B]). *Suppose 2 is linear. Let $G = \mathrm{GL}_{2^r}(q)$ where $r \geq 1$. Let $R = EZ$ be a subgroup of $G$ of symplectic-type, where $E \cong 2^{2r+1}_\eta$ and $Z$ is the Sylow 2-subgroup of the center of $G$. If we set $C = C_G(R)$ and $N = N_G(R)$, then $C_N(R) = Z(N) = Z(G)$ and $N/RC \cong \mathrm{Sp}_{2r}(2)$.*

By the two propositions, $\bar{\Gamma}_r.\mathrm{Sp}_{2r}(2)$ is embedded in $\mathrm{PGL}_{2^r}(q)$ for $q$ an even power of the prime $\ell$. Therefore, $F$ acts trivially on $N/C^\circ$ when $\Lambda$ is a nontoral elementary abelian 2-subgroup of the algebraic group $\mathrm{PGL}_n(K)$. Hence the conclusions of Sections 4.2 and 4.3 on the splittings and their local structure in the finite group $\mathrm{PGL}_n(q)$ hold.

This completes the discussion of the elementary abelian 2-subgroups of $\mathrm{PGL}_n(q)$. Next, we consider the elementary abelian 2-subgroups of $\mathrm{PGU}_n(q)$.

### 4.5. Nontoral 2-subgroups of $\mathrm{PGU}_n(q)$

In this section, we discuss the elementary abelian 2-subgroups of $\mathrm{PGU}_n(q)$. We consider the case where $q \equiv 3 \bmod 4$ in Section 4.5.1 and $q \equiv 1 \bmod 4$ in Section 4.5.2. If $q \equiv 3 \bmod 4$, then the $F$-action is trivial; otherwise, it is nontrivial. In each case, we discuss separately the cases where $\bar{A}$ has a connected and disconnected centralizer.

#### 4.5.1. $\mathrm{PGU}_n(q)$ where $q \equiv 3 \bmod 4$

Let $K$ be an algebraically closed field of characteristic $\ell$. Let $F$ be a Steinberg endomorphism of $\mathrm{PGL}_n(K)$ such that $\mathrm{PGL}_n(K)^F = \mathrm{PGU}_n(q)$ where $q$ is a power of $\ell$ and $q \equiv 3 \bmod 4$. Let $\Lambda$ be a nontoral elementary abelian 2-subgroup of $\mathrm{PGU}_n(q)$. Let $n = 2^r \times k$ where both $r$ and $k$ are positive integers. Let $E \leq \mathrm{PGL}_n(\mathbb{C})$ correspond to $\Lambda$ under the bijection of Proposition 3.1. So

$$E = \bar{\Gamma}_r \times \bar{A}$$

where $\bar{A}$ is either trivial or a representative of a conjugacy class of the toral elementary abelian 2-subgroups of $\mathrm{PGL}_k(\mathbb{C})$.

By Proposition 3.2, $\Lambda = 2^{2r} \times \Lambda_w$ where $\Lambda_w$ is either trivial or a toral elementary abelian 2-subgroup of $\mathrm{PGL}_k(K)$. We denote $N_{\mathrm{PGL}_n(K)}(\Lambda)$ by $N$ and $C_{\mathrm{PGL}_n(K)}(\Lambda)$ by $C$.

When $\bar{A}$ is trivial, $E = \bar{\Gamma}_r$, and $C_{\mathrm{PGL}_n(\mathbb{C})}(E) \cong \bar{\Gamma}_r \times \mathrm{PGL}_k(\mathbb{C})$ and $W_{\mathrm{PGL}_n(\mathbb{C})}(E) = \mathrm{Sp}_{2r}(2)$. By [2, 1F], for $r \geq 1$, the group $\bar{\Gamma}_r.\mathrm{Sp}_{2r}(2)$ is a subgroup of $\mathrm{PGU}_{2^r}(q)$ for $q \equiv 3 \bmod 4$ and is therefore centralized by $F$. Hence, when $\Lambda = 2^{2r}$, there are two $F$-classes of $N/C^\circ$ in $C/C^\circ$, one consisting of the identity and the other all the remaining elements. Correspondingly, there are two $\mathrm{PGU}_n(q)$-conjugacy classes of elementary abelian 2-subgroups with representatives $E_1$ and $E_2$. Their local structure

is the following:
$$\begin{aligned}
C_{\mathrm{PGU}_n(q)}(E_1) &\cong \bar{\Gamma}_r \times \mathrm{PGU}_k(q) \\
N_{\mathrm{PGU}_n(q)}(E_1) &\cong (\bar{\Gamma}_r \times \mathrm{PGU}_k(q)).\mathrm{Sp}_{2r}(2) \\
C_{\mathrm{PGU}_n(q)}(E_2) &\cong \bar{\Gamma}_r \times \mathrm{PGU}_k(q) \\
N_{\mathrm{PGU}_n(q)}(E_2) &\cong (\bar{\Gamma}_r \times \mathrm{PGU}_k(q)).(2^{2r-1}:\mathrm{Sp}_{2(r-1)}(2)).
\end{aligned}$$

Now we discuss the nontrivial $\bar{A}$ with connected centralizers. Note that $F$ acts trivially on $N/C^\circ$ when $q \equiv 3 \bmod 4$.

**Case 1**. $\bar{A}$ is nontrivial toral in $\mathrm{PGL}_k(\mathbb{C})$ with $C_{\mathrm{PGL}_k(\mathbb{C})}(\bar{A})$ connected.

In this case,
$$\begin{aligned}
C_{\mathrm{PGL}_n(\mathbb{C})}(E)/C_{\mathrm{PGL}_n(\mathbb{C})}(E)^\circ &\cong \bar{\Gamma}_r \\
W_{\mathrm{PGL}_n(\mathbb{C})}(E) &= 2^{\mathrm{rk}\bar{A} \times 2r} : (\mathrm{Sp}_{2r}(2) \times W_{\mathrm{PGL}_k(\mathbb{C})}(\bar{A})).
\end{aligned}$$

By construction of the generators of $U_1$ from Theorem 4.6, $C_{\mathrm{PGL}_n(\mathbb{C})}(E)/C_{\mathrm{PGL}_n(\mathbb{C})}(E)^\circ$ is centralized by $U_1$. Therefore, there are two $F$-classes of $N/C^\circ$ in $C/C^\circ$, one consisting of the identity and the other all the remaining elements of $2^{2r}$, the first component of $\Lambda$. Correspondingly, there are two $\mathrm{PGU}_n(q)$-conjugacy classes of elementary abelian 2-subgroups with representatives $E_1$ and $E_2$. Their local structure is the following:

$$\begin{aligned}
C_{\mathrm{PGU}_n(q)}(E_1) &\cong \bar{\Gamma}_r \times \left(C_{\mathrm{PGL}_k(K)}(\Lambda_w)\right)^F \\
C_{\mathrm{PGU}_n(q)}(E_2) &\cong \bar{\Gamma}_r \times \left(C_{\mathrm{PGL}_k(K)}(\Lambda_w)\right)^F \\
N_{\mathrm{PGU}_n(q)}(E_1) &\cong (\bar{\Gamma}_r \times \left(C_{\mathrm{PGL}_k(K)}(\Lambda_w)\right)^F).W_1 \\
N_{\mathrm{PGU}_n(q)}(E_2) &\cong (\bar{\Gamma}_r \times \left(C_{\mathrm{PGL}_k(K)}(\Lambda_w)\right)^F).W_2
\end{aligned}$$

where
$$\begin{aligned}
W_1 &= 2^{\mathrm{rk}\bar{A} \times 2r} : (\mathrm{Sp}_{2r}(2) \times W_{\mathrm{PGL}_k(\mathbb{C})}(\bar{A})) \\
W_2 &= 2^{\mathrm{rk}\bar{A} \times 2r} : ((2^{2r-1}:\mathrm{Sp}_{2(r-1)}(2)) \times W_{\mathrm{PGL}_k(\mathbb{C})}(\bar{A})).
\end{aligned}$$

Next, we discuss the nontrivial $\bar{A}$ with disconnected centralizers.

**Case 2**. $\bar{A}$ is toral in $\mathrm{PGL}_k(\mathbb{C})$ with $C_{\mathrm{PGL}_k(\mathbb{C})}(\bar{A})$ disconnected.

Let $\Lambda = 2^{2r} \times \Lambda_w$ be a representative of a conjugacy class of the nontoral elementary abelian 2-subgroups of $G := \mathrm{PGL}_n(K)$ where $\Lambda_w$ is a representative of a conjugacy class of the toral elementary abelian 2-subgroups of $\mathrm{PGL}_k(K)$ with $C_{\mathrm{PGL}_k(K)}(\Lambda_w)$ disconnected.

For each factorization $k = 2^{r_1} \times k_1$ where both $r_1$ and $k_1$ are positive integers, $\Lambda_w = D_{r_1} \times H$. Here $D_{r_1} = \langle \bar{A}_0^k, ..., \bar{A}_{r-1}^k \rangle$ where the superscript $k$ means elements in $\operatorname{PGL}_k(K)$, and $H$ is either trivial or a representative of a conjugacy class of the toral elementary abelian 2-subgroups of $\operatorname{PGL}_{k_1}(K)$ with $C_{\operatorname{PGL}_{k_1}(K)}(H)$ connected. Let $B_{r_1} = \langle \bar{B}_0^k, ..., \bar{B}_{r-1}^k \rangle$. Now

$$\begin{aligned} C_G(\Lambda) &\cong \bar{\Gamma}_r \times (C_{\operatorname{PGL}_k(K)}(\Lambda_w)^\circ . B_{r_1}) \\ N_G(\Lambda)/C_G(\Lambda) &\cong 2^{\operatorname{rk}\Lambda_w \times 2r} : (\operatorname{Sp}_{2r}(2) \times W_{\operatorname{PGL}_k(K)}(\Lambda_w)) \end{aligned}$$

where

$$W_{\operatorname{PGL}_k(K)}(\Lambda_w) \cong 2^{\operatorname{rk}H \times r_1} : (\operatorname{GL}_{r_1}(2) \times W_{\operatorname{PGL}_{k_1}(K)}(H)).$$

Now $F$ centralizes $N/C^\circ$ when $q \equiv 3 \mod 4$. Therefore, there are three $F$-classes of $N/C^\circ$ in $C/C^\circ$ when $\Lambda = 2^{2r} \times \Lambda_w$ where $\Lambda_w$ has a disconnected centralizer in $\operatorname{PGL}_k(K)$. One consists of the identity, another consists of the nontrivial elements of $2^{2r}$ and the third consists of all the remaining elements with representative $b$ (see Case 2 of Section 4.3.3). Correspondingly, there are three $\operatorname{PGU}_n(q)$-conjugacy classes of elementary abelian 2-subgroups with representatives $E_1$, $E_2$ and $E_3$. Their local structure is the following:

$$\begin{aligned} C_{\operatorname{PGU}_n(q)}(E_1) &\cong \bar{\Gamma}_r \times ((C_{\operatorname{PGL}_k(K)}(\Lambda_w)^\circ)^F . B_{r_1}) \\ C_{\operatorname{PGU}_n(q)}(E_2) &\cong \bar{\Gamma}_r \times ((C_{\operatorname{PGL}_k(K)}(\Lambda_w)^\circ)^F . B_{r_1}) \\ C_{\operatorname{PGU}_n(q)}(E_3) &\cong \bar{\Gamma}_r \times ((C_{\operatorname{PGL}_k(K)}(\Lambda_w)^\circ)^{bF} . B_{r_1}) \\ N_{\operatorname{PGU}_n(q)}(E_1) &\cong C_{\operatorname{PGU}_n(q)}(E_1) . (2^{\operatorname{rk}\Lambda_w \times 2r} : (\operatorname{Sp}_{2r}(2) \times W_{\operatorname{PGL}_k(K)}(\Lambda_w))) \\ N_{\operatorname{PGU}_n(q)}(E_2) &\cong C_{\operatorname{PGU}_n(q)}(E_2) . (2^{\operatorname{rk}\Lambda_w \times 2r} : ((2^{2r-1} : \operatorname{Sp}_{2(r-1)}(2)) \times W_{\operatorname{PGL}_k(K)}(\Lambda_w))) \\ N_{\operatorname{PGU}_n(q)}(E_3) &\cong C_{\operatorname{PGU}_n(q)}(E_3) . (2^{(\operatorname{rk}\Lambda_w - 1) \times 2r} : (\operatorname{Sp}_{2r}(2) \times W_3)) \end{aligned}$$

where

$$W_3 \cong 2^{\operatorname{rk}H \times r_1} : (\begin{pmatrix} 1 & * \\ 0 & \operatorname{GL}_{r_1-1}(2) \end{pmatrix} \times W_{\operatorname{PGL}_{k_1}(K)}(H)).$$

4.5.2. $\operatorname{PGU}_n(q)$ *where* $q \equiv 1 \mod 4$

Let $K$ be an algebraically closed field of characteristic $\ell$. Let $F$ be a Steinberg endomorphism of $\operatorname{PGL}_n(K)$ such that $\operatorname{PGL}_n(K)^F = \operatorname{PGU}_n(q)$ where $q$ is a power of $\ell$ and $q \equiv 1 \mod 4$. Let $\Lambda$ be a nontoral elementary abelian 2-subgroup of $\operatorname{PGU}_n(q)$. Let $n = 2^r \times k$ where both $r$ and $k$ are positive integers. Let $E \leq \operatorname{PGL}_n(\mathbb{C})$ correspond to $\Lambda$ under the bijection of Proposition 3.1. So

$$E = \bar{\Gamma}_r \times \bar{A}$$

where $\bar{A}$ is either trivial or a representative of a conjugacy class of the toral elementary abelian 2-subgroups of $\operatorname{PGL}_k(\mathbb{C})$.

By Proposition 3.2, $\Lambda = 2^{2r} \times \Lambda_w$ where $\Lambda_w$ is either trivial or a toral elementary abelian 2-subgroup of $\operatorname{PGL}_k(K)$. We denote $N_{\operatorname{PGL}_n(K)}(\Lambda)$ by $N$ and $C_{\operatorname{PGL}_n(K)}(\Lambda)$ by

$C$.

The following proposition is used when we descend to the finite groups. The proof is similar to that of Theorem 4.17.

**Proposition 4.20.** *Let $H = \mathrm{PGL}_{2^r}(K)$ where $r \geq 1$. Let $F$ be a Steinberg endomorphism such that $H^F = \mathrm{PGU}_{2^r}(q)$ where $q \equiv 1 \bmod 4$. Let $E$ be a nontoral elementary abelian subgroup of $H$ of order $2^{2r}$. There exist two $F$-classes of $N_H(E)/C_H(E)^\circ$ in $C_H(E)/C_H(E)^\circ$. Correspondingly, there are two $H^F$-conjugacy classes of elementary abelian 2-subgroups of order $2^{2r}$ with representatives $E_1$ and $E_2$. Their local structure is the following:*

$$C_{H^F}(E_1) \cong 2^{2r} \text{ and } N_{H^F}(E_1) \cong 2^{2r}.\mathrm{SO}_{2r}^+(2)$$

$$C_{H^F}(E_2) \cong 2^{2r} \text{ and } N_{H^F}(E_2) \cong 2^{2r}.\mathrm{SO}_{2r}^-(2).$$

Next, we descend to the finite groups.

When $\bar{A}$ is trivial, $E = \bar{\Gamma}_r$, and $C_{\mathrm{PGL}_n(\mathbb{C})}(E) \cong \bar{\Gamma}_r \times \mathrm{PGL}_k(\mathbb{C})$ and $W_{\mathrm{PGL}_n(\mathbb{C})}(E) = \mathrm{Sp}_{2r}(2)$. Hence, by Proposition 4.20, when $\Lambda = 2^{2r}$, there are two $F$-classes of $N/C^\circ$ in $C/C^\circ$, and correspondingly, there are two $\mathrm{PGU}_n(q)$-conjugacy classes of elementary abelian 2-subgroups with representatives $E_1$ and $E_2$. Their local structure is the following:

$$\begin{aligned}
C_{\mathrm{PGU}_n(q)}(E_1) &\cong \bar{\Gamma}_r \times \mathrm{PGU}_k(q) \\
N_{\mathrm{PGU}_n(q)}(E_1) &\cong (\bar{\Gamma}_r \times \mathrm{PGU}_k(q)).\mathrm{SO}_{2r}^+(2) \\
C_{\mathrm{PGU}_n(q)}(E_2) &\cong \bar{\Gamma}_r \times \mathrm{PGU}_k(q) \\
N_{\mathrm{PGU}_n(q)}(E_2) &\cong (\bar{\Gamma}_r \times \mathrm{PGU}_k(q)).\mathrm{SO}_{2r}^-(2).
\end{aligned}$$

Now we focus on the nontrivial $\bar{A}$ and consider two cases.

**Case 1**. $\bar{A}$ is nontrivial toral in $\mathrm{PGL}_k(\mathbb{C})$ with $C_{\mathrm{PGL}_k(\mathbb{C})}(\bar{A})$ connected.

A representative of a conjugacy class of the nontoral elementary abelian 2-subgroups of $\mathrm{PGL}_n(\mathbb{C})$ is

$$E = \bar{\Gamma}_r \times \bar{A}$$

where $\bar{A}$ is a representative of a conjugacy class of the nontrivial toral elementary abelian 2-subgroups of $\mathrm{PGL}_k(\mathbb{C})$ and $C_{\mathrm{PGL}_k(\mathbb{C})}(\bar{A})$ is connected.

Since $C_{\mathrm{PGL}_k(\mathbb{C})}(\bar{A})$ is connected, $C_{\mathrm{PGL}_n(\mathbb{C})}(E)/C_{\mathrm{PGL}_n(\mathbb{C})}(E)^\circ \cong \bar{\Gamma}_r$. By construction of the generators of $U_1$ from Theorem 4.6, $C_{\mathrm{PGL}_n(\mathbb{C})}(E)/C_{\mathrm{PGL}_n(\mathbb{C})}(E)^\circ$ is centralized by $U_1$. By Proposition 4.20, when $\Lambda = 2^{2r} \times \Lambda_w$ where $\Lambda_w$ has a connected centralizer in $\mathrm{PGL}_k(K)$, there are two $F$-classes of $N/C^\circ$ in $C/C^\circ$. Correspondingly, there are two $\mathrm{PGU}_n(q)$-conjugacy classes of elementary abelian 2-subgroups with representatives

$E_1$ and $E_2$. Their local structure is the following:

$$\begin{aligned}
C_{\mathrm{PGU}_n(q)}(E_1) &\cong \bar{\Gamma}_r \times (C_{\mathrm{PGL}_k(K)}(\Lambda_w))^F \\
C_{\mathrm{PGU}_n(q)}(E_2) &\cong \bar{\Gamma}_r \times (C_{\mathrm{PGL}_k(K)}(\Lambda_w))^F \\
N_{\mathrm{PGU}_n(q)}(E_1) &\cong (\bar{\Gamma}_r \times (C_{\mathrm{PGL}_k(K)}(\Lambda_w))^F).W_1 \\
N_{\mathrm{PGU}_n(q)}(E_2) &\cong (\bar{\Gamma}_r \times (C_{\mathrm{PGL}_k(K)}(\Lambda_w))^F).W_2
\end{aligned}$$

where

$$\begin{aligned}
W_1 &= 2^{\mathrm{rk}\bar{A}\times 2r} : (\mathrm{SO}_{2r}^+(2) \times W_{\mathrm{PGL}_k(\mathbb{C})}(\bar{A})) \\
W_2 &= 2^{\mathrm{rk}\bar{A}\times 2r} : (\mathrm{SO}_{2r}^-(2) \times W_{\mathrm{PGL}_k(\mathbb{C})}(\bar{A})).
\end{aligned}$$

**Case 2**. $\bar{A}$ is toral in $\mathrm{PGL}_k(\mathbb{C})$ with $C_{\mathrm{PGL}_k(\mathbb{C})}(\bar{A})$ disconnected.

Let $\Lambda = 2^{2r} \times \Lambda_w$ be a representative of a conjugacy class of the nontoral elementary abelian 2-subgroups of $G := \mathrm{PGL}_n(K)$ where $\Lambda_w$ is a representative of a conjugacy class of toral elementary abelian 2-subgroups of $\mathrm{PGL}_k(K)$ with $C_{\mathrm{PGL}_k(K)}(\Lambda_w)$ disconnected.

For each factorization $k = 2^{r_1} \times k_1$ where both $r_1$ and $k_1$ are positive integers, $\Lambda_w = D_{r_1} \times H$. Here $D_{r_1} = \langle \bar{A}_0^k, ..., \bar{A}_{r-1}^k \rangle$ where the superscript $k$ means elements in $\mathrm{PGL}_k(K)$, and $H$ is either trivial or a representative of a conjugacy class of the toral elementary abelian 2-subgroups of $\mathrm{PGL}_{k_1}(K)$ with $C_{\mathrm{PGL}_{k_1}(K)}(H)$ connected. Let $B_{r_1} = \langle \bar{B}_0^k, ..., \bar{B}_{r-1}^k \rangle$. Now

$$\begin{aligned}
C_G(\Lambda) &\cong \bar{\Gamma}_r \times (C_{\mathrm{PGL}_k(K)}(\Lambda_w)^\circ.B_{r_1}) \\
N_G(\Lambda)/C_G(\Lambda) &\cong 2^{\mathrm{rk}\Lambda_w \times 2r} : (\mathrm{Sp}_{2r}(2) \times W_{\mathrm{PGL}_k(K)}(\Lambda_w))
\end{aligned}$$

where

$$W_{\mathrm{PGL}_k(K)}(\Lambda_w) \cong 2^{\mathrm{rk}H \times r_1} : (\mathrm{GL}_{r_1}(2) \times W_{\mathrm{PGL}_{k_1}(K)}(H)).$$

So

$$C_G(\Lambda)/C_G(\Lambda)^\circ \cong \bar{\Gamma}_r \times B_{r_1}.$$

When $\Lambda = 2^{2r} \times \Lambda_w$ where $\Lambda_w$ has a disconnected centralizer in $\mathrm{PGL}_k(K)$, there are three $F$-classes of $N/C^\circ$ in $C/C^\circ$. One of the three classes is $Or_1' \cup Or_2'$ with representative $b'$ (see Case 2 of Section 4.4.) Correspondingly, there are three $\mathrm{PGU}_n(q)$-conjugacy classes of elementary abelian 2-subgroups with representatives $E_1$, $E_2$ and

$E_3$. Their local structure is the following:

$$\begin{align}
C_{\mathrm{PGU}_n(q)}(E_1) &\cong \bar{\Gamma}_r \times ((C_{\mathrm{PGL}_k(K)}(\Lambda_w)^\circ)^F.B_{r_1}) \\
C_{\mathrm{PGU}_n(q)}(E_2) &\cong \bar{\Gamma}_r \times ((C_{\mathrm{PGL}_k(K)}(\Lambda_w)^\circ)^F.B_{r_1}) \\
C_{\mathrm{PGU}_n(q)}(E_3) &\cong \bar{\Gamma}_r \times ((C_{\mathrm{PGL}_k(K)}(\Lambda_w)^\circ)^{b'F}.B_{r_1}) \\
N_{\mathrm{PGU}_n(q)}(E_1) &\cong C_{\mathrm{PGU}_n(q)}(E_1).(2^{\mathrm{rk}\Lambda_w \times 2r} : (\mathrm{SO}_{2r}^+(2) \times W_{\mathrm{PGL}_k(K)}(\Lambda_w))) \\
N_{\mathrm{PGU}_n(q)}(E_2) &\cong C_{\mathrm{PGU}_n(q)}(E_2).(2^{\mathrm{rk}\Lambda_w \times 2r} : (\mathrm{SO}_{2r}^-(2) \times W_{\mathrm{PGL}_k(K)}(\Lambda_w))) \\
N_{\mathrm{PGU}_n(q)}(E_3) &\cong C_{\mathrm{PGU}_n(q)}(E_3).(2^{(\mathrm{rk}\Lambda_w-1) \times 2r} : (\mathrm{Sp}_{2r}(2) \times W_3))
\end{align}$$

where

$$W_3 \cong 2^{\mathrm{rk}H \times r_1} : (\begin{pmatrix} 1 & * \\ 0 & \mathrm{GL}_{r_1-1}(2) \end{pmatrix} \times W_{\mathrm{PGL}_{k_1}(K)}(H)).$$

## 4.6. Nontoral elementary abelian $p$-subgroups of $\mathrm{PGL}_n(q)$ where $p$ is odd

In Theorem 4.23 and Proposition 4.25, we first examine the local structure of nontoral elementary abelian $p$-subgroups of the algebraic groups. Next, we obtain the structure and local structure of toral elementary abelian $p$-subgroups with disconnected centralizers in Theorems 4.27, 4.28 and 4.32. Lastly, we descend to the finite groups and discuss separately the cases where $\bar{A}$ has a connected and disconnected centralizer.

Recall the definition of symplectic-type $p$-groups from Section 4.4. Symplectic-type $p$-groups of exponent $p \cdot \gcd(2, p)$ are listed in [17, Table 4.6 A]. Note that representations $\phi_1, \phi_2$ of a finite group $X$ are *quasiequivalent* if there is an automorphism $\theta$ of $X$ such that $\phi_1$ is equivalent to $\theta \phi_2$.

**Proposition 4.21** ([17, Proposition 4.6.3]). *Assume that $R$ is a symplectic-type $p$-group of exponent $p \cdot \gcd(2, p)$ and order $p^{1+2m}$ where $m \geq 1$. Let $\ell$ be a prime such that $p \neq \ell$.*
  (i) *$R$ has precisely $|Z(R)| - 1$ inequivalent faithful absolutely irreducible representations over an algebraically closed field of characteristic $\ell$. Denote them by $\rho_1, ..., \rho_k$, where $k = |Z(R)| - 1$.*
  (ii) *The $\rho_i$ are quasiequivalent, they have degree $p^m$, and the smallest field over which they can be realized is $\mathrm{GF}(\ell^e)$, where $e$ is the smallest integer for which $\ell^e \equiv 1 \mod |Z(R)|$.*
  (iii) *If $i \neq j$, then $\rho_i$ and $\rho_j$ differ on $Z(R)$.*

**Lemma 4.22.** *Let $G = \mathrm{GL}_n(K)$ where $K$ is an algebraically closed field of characteristic $\ell$, let $F$ be a Steinberg endomorphism of $G$ and let $Y$ be an abelian $p$-subgroup of $G^F$ where $p \neq \ell$. The group $Y$ is toral in $G^F$.*

*Proof.* By [11, Proposition 2.1], $H := C_G(Y)^\circ$ is a Levi subgroup of $G$. Hence $Z(H)$ is contained in the intersection of all the maximal tori of $H$. Since $Y \leq Z(H)$, we deduce that $Y$ is toral. $\square$

Let $p$ be an odd prime and let $q$ be a power of a prime $\ell$ where $p \neq \ell$. Let $K$ be an algebraically closed field of characteristic $\ell$ and let $F$ be a Steinberg endomorphism of $\mathrm{PGL}_n(K)$ such that $\mathrm{PGL}_n(K)^F = \mathrm{PGL}_n(q)$. The following theorem is useful in the investigation of the $F$-action.

**Theorem 4.23.** *Let $p$ be an odd prime and $m \geq 1$. If there exists a nontoral elementary abelian subgroup $\Lambda$ of order $p^{2m}$ of $\mathrm{PGL}_n(q)$ with*
$$C_{\mathrm{PGL}_n(q)}(\Lambda) \cong \Lambda \times \mathrm{PGL}_k(q)$$
*where $1 \leq k < n$ and*
$$N_{\mathrm{PGL}_n(q)}(\Lambda) \cong C_{\mathrm{PGL}_n(q)}(\Lambda).\mathrm{Sp}_{2m}(p),$$
*then $p^m | n$ and $q \equiv 1 \bmod p$. Conversely, if $p^m | n$ and $q \equiv 1 \bmod p$, then there exists a nontoral elementary abelian subgroup $\Lambda$ of order $p^{2m}$ of $\mathrm{PGL}_n(q)$ with*
$$C_{\mathrm{PGL}_n(q)}(\Lambda) \cong \Lambda \times \mathrm{PGL}_k(q)$$
*where $1 \leq k < n$ and*
$$N_{\mathrm{PGL}_n(q)}(\Lambda) \cong C_{\mathrm{PGL}_n(q)}(\Lambda).\mathrm{Sp}_{2m}(p).$$

*Proof.* If $p^m | n$ and $q \equiv 1 \bmod p$, then by Proposition 4.21
$$p_+^{1+2m}.\mathrm{Sp}_{2m}(p) \leq \mathrm{GL}_n(q)$$
and
$$Z(p_+^{1+2m}.\mathrm{Sp}_{2m}(p)) \leq Z(\mathrm{GL}_n(q)).$$
If we denote by $\Lambda$ the image of $p_+^{1+2m}$ in $\mathrm{PGL}_n(q)$, then $\Lambda.\mathrm{Sp}_{2m}(p) \leq \mathrm{PGL}_n(q)$. Thus,
$$\Lambda \cong p^{2m} \leq \mathrm{PGL}_n(q) = \mathrm{PGL}_n(K)^F \leq \mathrm{PGL}_n(K)$$
and $\Lambda$ is nontoral in $\mathrm{PGL}_n(K)$ with
$$N_{\mathrm{PGL}_n(K)}(\Lambda)/C_{\mathrm{PGL}_n(K)}(\Lambda) \cong \mathrm{Sp}_{2m}(p)$$
since $\mathrm{Sp}_{2m}(p)$ is the outer automorphism group of $\Lambda$.

By the bijection of Proposition 3.1, let $S \cong p^{2m}$ be the nontoral subgroup of $\mathrm{PGL}_n(\mathbb{C})$ corresponding to $\Lambda$. By Proposition 3.2,
$$N_{\mathrm{PGL}_n(\mathbb{C})}(S)/C_{\mathrm{PGL}_n(\mathbb{C})}(S) \cong \mathrm{Sp}_{2m}(p).$$

By Theorem 4.4, if $E$ is a nontoral elementary abelian $p$-subgroup of $\mathrm{PGL}_n(\mathbb{C})$, then
$$E = \bar{\Gamma}_r \times \bar{A}$$
where $\bar{A}$ is either trivial or a representative of a conjugacy class of the toral elementary abelian $p$-subgroups of $\mathrm{PGL}_k(\mathbb{C})$. Also,
$$C_{\mathrm{PGL}_n(\mathbb{C})}(E) \cong \bar{\Gamma}_r \times C_{\mathrm{PGL}_k(\mathbb{C})}(\bar{A})$$

and
$$W_{\mathrm{PGL}_n(\mathbb{C})}(E) = p^{\mathrm{rk}\bar{A}\times 2r} : (\mathrm{Sp}_{2r}(p) \times W_{\mathrm{PGL}_k(\mathbb{C})}(\bar{A})).$$
Therefore, we deduce that $\bar{A}$ is trivial and $r = m$. Hence
$$C_{\mathrm{PGL}_n(\mathbb{C})}(S) \cong S \times \mathrm{PGL}_k(\mathbb{C})$$
and, by Proposition 3.2,
$$C_{\mathrm{PGL}_n(K)}(\Lambda) \cong \Lambda \times \mathrm{PGL}_k(K)$$
and
$$C_{\mathrm{PGL}_n(q)}(\Lambda) \cong \Lambda \times \mathrm{PGL}_k(q).$$
For the other direction, suppose there exists a nontoral elementary abelian subgroup $\Lambda$ of order $p^{2m}$ in $\mathrm{PGL}_n(q)$ with
$$C_{\mathrm{PGL}_n(q)}(\Lambda) \cong \Lambda \times \mathrm{PGL}_k(q)$$
where $k < n$ and
$$N_{\mathrm{PGL}_n(q)}(\Lambda) \cong C_{\mathrm{PGL}_n(q)}(\Lambda).\mathrm{Sp}_{2m}(p).$$
Since
$$\Lambda \leq \mathrm{PGL}_n(q) = \mathrm{PGL}_n(K)^F,$$
we deduce that $C_{\mathrm{PGL}_n(K)}(\Lambda) \cong \Lambda \times \mathrm{PGL}_k(K)$, and $N_{\mathrm{PGL}_n(K)}(\Lambda) \cong C_{\mathrm{PGL}_n(K)}(\Lambda).\mathrm{Sp}_{2m}(p)$ and $\Lambda$ is nontoral in $\mathrm{PGL}_n(K)$.

By the bijection of Proposition 3.1, let $S \cong p^{2m}$ be the subgroup of $\mathrm{PGL}_n(\mathbb{C})$ corresponding to $\Lambda$. By Proposition 3.2,
$$C_{\mathrm{PGL}_n(\mathbb{C})}(S) \cong S \times \mathrm{PGL}_k(\mathbb{C})$$
and
$$N_{\mathrm{PGL}_n(\mathbb{C})}(S)/C_{\mathrm{PGL}_n(\mathbb{C})}(S) \cong \mathrm{Sp}_{2m}(p).$$
By Theorem 4.4, if $E$ is a nontoral elementary abelian $p$-subgroup of $\mathrm{PGL}_n(\mathbb{C})$, then
$$E = \bar{\Gamma}_r \times \bar{A}$$
where $\bar{A}$ is either trivial or a representative of a conjugacy class of the toral elementary abelian $p$-subgroups of $\mathrm{PGL}_k(\mathbb{C})$. Also,
$$C_{\mathrm{PGL}_n(\mathbb{C})}(E) \cong \bar{\Gamma}_r \times C_{\mathrm{PGL}_k(\mathbb{C})}(\bar{A})$$
and
$$W_{\mathrm{PGL}_n(\mathbb{C})}(E) = p^{\mathrm{rk}\bar{A}\times 2r} : (\mathrm{Sp}_{2r}(p) \times W_{\mathrm{PGL}_k(\mathbb{C})}(\bar{A})).$$
Therefore, we deduce that $S$ is nontoral, $r = m$, $\bar{A}$ is trivial and $p^m | n$.

Finally, suppose $q \not\equiv 1 \bmod p$. Let $Y_1$ be a Sylow $p$-subgroup of $\mathrm{PGL}_n(q)$ and let $Y_2$ be a Sylow $p$-subgroup of $\mathrm{GL}_n(q)$. Since $Z(\mathrm{GL}_n(q)) \cong q-1$, we deduce that $Y_1 \cong Y_2$. By Lemma 4.22, all the abelian $p$-subgroups of $Y_2$ are toral; so are those of $Y_1$. Thus, we have a contradiction. $\square$

**Example 4.24.** *The local structure restrictions are needed. Since $5_+^{1+2} \leq \mathrm{GL}_5(11)$ and $C_5 \times C_5 \leq \mathrm{GL}_3(11)$, we can embed $R := 5_+^{1+2} \times C_5 \times C_5$ in $\mathrm{GL}_{15}(11)$ and hence embed $R/Z(R) \cong 5^4$ in $\mathrm{PGL}_{15}(11)$. Computation shows $N_{\mathrm{PGL}_{15}(11)}(E)/C_{\mathrm{PGL}_{15}(11)}(E) \not\cong \mathrm{Sp}_4(5)$. But $5^2 \nmid 15$.*

4.6.1. *The $F$-action on the normalizer quotients of nontoral $p$-subgroups.*

Let $\Lambda$ be a nontoral elementary abelian $p$-subgroup of $\mathrm{PGL}_n(K)$ where $n = p^r \times k$ and both $r$ and $k$ are positive integers. Let $\mathrm{PGL}_n(K)^F = \mathrm{PGL}_n(q)$ where $q \equiv 1 \bmod p$ and $\Lambda \leq \mathrm{PGL}_n(K)^F$. We denote $N_{\mathrm{PGL}_n(K)}(\Lambda)$ by $N$ and $C_{\mathrm{PGL}_n(K)}(\Lambda)$ by $C$.

Let $E$ be the nontoral elementary abelian $p$-subgroup of $\mathrm{PGL}_n(\mathbb{C})$ corresponding to $\Lambda$ under the bijection of Proposition 3.1. By Theorem 4.4,

$$E = \bar{\Gamma}_r \times \bar{A}$$

where $\bar{A}$ is either trivial or a representative of a conjugacy class of the toral elementary abelian $p$-subgroups of $\mathrm{PGL}_k(\mathbb{C})$. By Proposition 3.2, $\Lambda = p^{2r} \times \Lambda_w$ where $\Lambda_w$ is either trivial or a toral elementary abelian $p$-subgroup of $\mathrm{PGL}_k(K)$. Now

$$C_{\mathrm{PGL}_n(\mathbb{C})}(E) \cong \bar{\Gamma}_r \times C_{\mathrm{PGL}_k(\mathbb{C})}(\bar{A})$$

and

$$W_{\mathrm{PGL}_n(\mathbb{C})}(E) = p^{\mathrm{rk}\bar{A} \times 2r} : (\mathrm{Sp}_{2r}(p) \times W_{\mathrm{PGL}_k(\mathbb{C})}(\bar{A})).$$

Let $J$ be a maximal toral elementary abelian $p$-subgroup of $\mathrm{PGL}_k(\mathbb{C})$. The rank of $J$ is $m := k - 1$. Let $Q = \bar{\Gamma}_r \times J$. Now

$$C_{\mathrm{PGL}_n(\mathbb{C})}(Q) \cong \bar{\Gamma}_r \times T_m$$

and

$$W_{\mathrm{PGL}_n(\mathbb{C})}(Q) = p^{m \times 2r} : (\mathrm{Sp}_{2r}(p) \times S_{m+1}).$$

We denote by $R$ the subgroup $p^{m \times 2r}$ of $W_{\mathrm{PGL}_n(\mathbb{C})}(Q)$. Let $A \in \mathrm{Sp}_{2r}(p)$, $B \in S_{m+1}$ and $u \in R$. Now

$$\begin{pmatrix} A & 0 \\ 0 & B \end{pmatrix} \begin{pmatrix} \mathrm{Id} & 0 \\ u & \mathrm{Id} \end{pmatrix} \begin{pmatrix} A^{-1} & 0 \\ 0 & B^{-1} \end{pmatrix} = \begin{pmatrix} \mathrm{Id} & 0 \\ BuA^{-1} & \mathrm{Id} \end{pmatrix}$$

implies that the generators of $R$ are in one orbit under the action of $\mathrm{Sp}_{2r}(p) \times S_{m+1}$. Let $\alpha$ be a primitive $p$-th root of unity in $\mathbb{C}$. Observe that $R$ has a generator

$$\left[ \left( \Delta\begin{pmatrix} 1 & & & & \\ & \alpha & & & \\ & & \alpha^2 & & \\ & & & \ddots & \\ & & & & \alpha^{p-1} \end{pmatrix}\right)_{p^{r-1}} \quad \Delta(I_{p^r})_m \right] C_{\mathrm{PGL}_n(\mathbb{C})}(Q).$$

Since all of the above nonzero entries are $1, \alpha, \alpha^2, ...,$ or $\alpha^{p-1}$, similarly to the case where $p = 2$, by Theorem 4.23, the Steinberg endomorphism $F$ centralizes $N/C^\circ$ when $\Lambda$ is a maximal nontoral elementary abelian $p$-subgroup of $\mathrm{PGL}_n(K)$, and

$q \equiv 1 \mod p$.

As in Section 4.2, we construct $2mr$ elements $x_{ij}, y_{ij} \in \mathrm{PGL}_n(\mathbb{C})$ for $1 \leq i \leq m$ and $1 \leq j \leq r$ where

$$x_{ij} = \left[\begin{pmatrix} \Delta(I_{p^r})_{i-1} & & \\ & \Delta\left(\begin{pmatrix} & I_{p^{j-1}} & & & \\ & & I_{p^{j-1}} & & \\ & & & \ddots & \\ & & & & I_{p^{j-1}} \\ I_{p^{j-1}} & & & & \end{pmatrix}\right)_{p^{r-j}} & \\ & & \Delta(I_{p^r})_{m-(i-1)} \end{pmatrix}\right]$$

$$y_{ij} = \left[\begin{pmatrix} \Delta(I_{p^r})_{i-1} & & \\ & \Delta\left(\begin{pmatrix} I_{p^{j-1}} & & & & \\ & \alpha I_{p^{j-1}} & & & \\ & & \alpha^2 I_{p^{j-1}} & & \\ & & & \ddots & \\ & & & & \alpha^{p-1} I_{p^{j-1}} \end{pmatrix}\right)_{p^{r-j}} & \\ & & \Delta(I_{p^r})_{m-(i-1)} \end{pmatrix}\right].$$

In particular,

$$\begin{pmatrix} & I_{p^{j-1}} & & & \\ & & I_{p^{j-1}} & & \\ & & & \ddots & \\ & & & & I_{p^{j-1}} \\ I_{p^{j-1}} & & & & \end{pmatrix} \text{ and } \begin{pmatrix} I_{p^{j-1}} & & & & \\ & \alpha I_{p^{j-1}} & & & \\ & & \alpha^2 I_{p^{j-1}} & & \\ & & & \ddots & \\ & & & & \alpha^{p-1} I_{p^{j-1}} \end{pmatrix}$$

are $p \times p$ block matrices.

In Section 4.1, we defined $\bar{\Gamma}_r = \langle \bar{A}_s, \bar{B}_s \mid 0 \leq s \leq r-1 \rangle$. Now $J$ has generators $g_i = [M_i]$ where $i = 1, ..., m$ and $M_i$ is a $k \times k$ diagonal matrix with $\alpha$ in the $i$-th position and $1$ in the remaining positions. Hence the maximal nontoral subgroup $Q = \bar{\Gamma}_r \times J$ of $\mathrm{PGL}_n(\mathbb{C})$ has generators $\bar{A}_s$, $\bar{B}_s$ and $g_i \otimes I_{p^r}$ for $0 \leq s \leq r-1$ and $i = 1, ..., m$.

Let

$$u = \begin{pmatrix} & 1 & & & \\ & & 1 & & \\ & & & \ddots & \\ & & & & 1 \\ 1 & & & & \end{pmatrix}.$$

Observe

$$u \begin{pmatrix} 1 & & & & \\ & \alpha & & & \\ & & \alpha^2 & & \\ & & & \ddots & \\ & & & & \alpha^{p-1} \end{pmatrix} u^{-1} = \alpha \begin{pmatrix} 1 & & & & \\ & \alpha & & & \\ & & \alpha^2 & & \\ & & & \ddots & \\ & & & & \alpha^{p-1} \end{pmatrix}$$

and

$$\begin{pmatrix} 1 & & & \\ & \alpha & & \\ & & \alpha^2 & \\ & & & \ddots & \\ & & & & \alpha^{p-1} \end{pmatrix} u \begin{pmatrix} 1 & & & \\ & \alpha^{p-1} & & \\ & & \alpha^{p-2} & \\ & & & \ddots & \\ & & & & \alpha \end{pmatrix} = \alpha^{p-1} u.$$

Hence

$$x_{ij}\bar{A}_s {x_{ij}}^{-1} = \begin{cases} \bar{A}_s & s+1 \neq j \\ \bar{A}_s(g_i \otimes I_{p^r}) & s+1 = j \end{cases}$$

$$x_{ij}\bar{B}_s {x_{ij}}^{-1} = \bar{B}_s \text{ for all } s$$

$$y_{ij}\bar{A}_s {y_{ij}}^{-1} = \bar{A}_s \text{ for all } s$$

$$y_{ij}\bar{B}_s {y_{ij}}^{-1} = \begin{cases} \bar{B}_s & s+1 \neq j \\ \bar{B}_s(g_i^{-1} \otimes I_{p^r}) & s+1 = j \end{cases}$$

$x_{ij}$ and $y_{ij}$ fix all the $g_i \otimes I_{p^r}$ for $i = 1, ..., m$.

Let
$$N_1 = \langle x_{ij}, y_{ij} \mid 1 \leq i \leq m \text{ and } 1 \leq j \leq r \rangle.$$

If $n \in N_1$ is nontrivial, then $n \in N_{\mathrm{PGL}_n(\mathbb{C})}(Q) \setminus C_{\mathrm{PGL}_n(\mathbb{C})}(Q)$. By the action of $N_1$ on $Q$ and the structure of $W_{\mathrm{PGL}_n(\mathbb{C})}(Q)$, we can identify $N_1 C_{\mathrm{PGL}_n(\mathbb{C})}(Q)/C_{\mathrm{PGL}_n(\mathbb{C})}(Q)$ with $R$, the subgroup of $W_{\mathrm{PGL}_n(\mathbb{C})}(Q)$. Thus, $N_1 T_m / T_m$ centralizes $\bar{\Gamma}_r$.

Let $E = \bar{\Gamma}_r \times L$ be a nontoral elementary abelian $p$-subgroup of $\mathrm{PGL}_n(\mathbb{C})$ where $L < J$. Now
$$C_{\mathrm{PGL}_n(\mathbb{C})}(E) \cong \bar{\Gamma}_r \times C_{\mathrm{PGL}_k(\mathbb{C})}(L)$$
and
$$W_{\mathrm{PGL}_n(\mathbb{C})}(E) = p^{\mathrm{rk}L \times 2r} : (\mathrm{Sp}_{2r}(p) \times W_{\mathrm{PGL}_k(\mathbb{C})}(L)).$$
We denote by $U_w$ the subgroup $p^{\mathrm{rk}L \times 2r}$ of $W_{\mathrm{PGL}_n(\mathbb{C})}(E)$.

The following proposition shows that the $F$-action is trivial.

**Proposition 4.25.** *Let $G = \mathrm{PGL}_n(\mathbb{C})$ where $n = p^r \times k$, $p$ is an odd prime and both $r$ and $k$ are positive integers. Let $Q = \bar{\Gamma}_r \times J$ be a maximal nontoral elementary abelian $p$-subgroup of $G$ where $J$ is a maximal toral elementary abelian $p$-subgroup of $\mathrm{PGL}_k(\mathbb{C})$ of rank $m$ where $m = k-1$. Let $E = \bar{\Gamma}_r \times L$ be a nontoral elementary abelian $p$-subgroup of $G$ where $L < J$. There exists $U_1 \leq N_1$ such that $u \in N_G(E) \setminus C_G(E)$ for each nontrivial $u \in U_1$ and $U_1 C_G(E)/C_G(E)$ is the subgroup $U_w$ of $W_G(E)$.*

*Proof.* The proof is similar to that for Theorem 4.6. □

By the above proposition, as in the maximal nontoral case, $N/C^\circ$ is centralized by $F$ when $\Lambda$ is a nonmaximal nontoral elementary abelian $p$-subgroup of $\mathrm{PGL}_k(K)$ and $q \equiv 1 \bmod p$.

### 4.6.2. *The (local) structure of toral p-subgroups with disconnected centralizers.*

Let $\Lambda = p^{2r} \times \Lambda_w$ where $\Lambda_w$ is toral and has a disconnected centralizer in $\mathrm{PGL}_k(K)$. We give the structure of such $\Lambda_w$ in Theorem 4.27. The local structure of $\Lambda_w$ is recorded in Theorem 4.28 and Theorem 4.32. First, we prove the following result which is needed for the structure of such $\Lambda_w$.

**Theorem 4.26.** *Let $G = \mathrm{PGL}_n(K)$ where $n = p \times k$, $p$ is an odd prime and $k \geq 1$. Let $\alpha \in K$ be a primitive p-th root of unity. There is a unique conjugacy class of toral elements of order $p$ in $G$ with disconnected centralizers. The class has representative*

$$e_{n/p} = \left[\begin{pmatrix} I_k & & & & \\ & \alpha I_k & & & \\ & & \alpha^2 I_k & & \\ & & & \ddots & \\ & & & & \alpha^{p-1} I_k \end{pmatrix}\right]$$

*and $C_G\left(e_{n/p}\right) = C_G(e_{n/p})^\circ . \langle f_{n/p} \rangle$ where*

$$f_{n/p} = \left[\begin{pmatrix} & I_k & & & \\ & & I_k & & \\ & & & \ddots & \\ & & & & I_k \\ I_k & & & & \end{pmatrix}\right].$$

*Proof.* Without loss of generality, let the maximal torus $T$ of $G$ be the image in $G$ of the group of diagonal matrices in $\mathrm{GL}_n(K)$.

Let $e \in T$ be a representative of a conjugacy class of the toral elements of order $p$ in $G$. Thus, $e$ is diagonal with entries $1, \alpha, \alpha^2, ..., \alpha^{p-2}$ or $\alpha^{p-1}$. The Weyl group $W = N_G(T)/T \cong S_n$. Suppose $\alpha^{i_1}, \alpha^{i_2}, ..., \alpha^{i_j}$ are the distinct powers of $\alpha$ appearing as entries on the diagonal of $e$ where $1 \leq j \leq p$ and $0 \leq i_1 < i_2 < ... < i_j \leq p-1$.

Since $C_G(e)$ is not connected, by Theorem 4.7, $\alpha^{i_1}, \alpha^{i_2}, ..., \alpha^{i_j}$ must have the same multiplicity.

Let $\sigma$ be a $j$-cycle and let $s_0 = \sigma^0(1) = 1, s_1 = \sigma(1), s_2 = \sigma^2(1), ..., s_{j-1} = \sigma^{j-1}(1)$. Thus,
$$\sigma = (s_0 \; s_1 \; s_2 \; ... \; s_{j-1}).$$
Now $\sigma \in W$. Let $\sigma$ centralize $e$. Hence $\sigma$ sends $(\alpha^{i_1}, \alpha^{i_2}, ..., \alpha^{i_j})$ to $\alpha^t(\alpha^{i_1}, \alpha^{i_2}, ..., \alpha^{i_j})$ for some integer $t$ where $1 \leq t < p$. We deduce that
$$\alpha^{i_{s_0}}/\alpha^{i_{s_1}} = \alpha^{i_{s_1}}/\alpha^{i_{s_2}} = ... = \alpha^{i_{s_{j-1}}}/\alpha^{i_{s_0}}$$
and $\alpha^{i_{s_1}}/\alpha^{i_{s_0}} = \alpha^t$. Thus,
$$1 = (\alpha^t)^j.$$
Since $\alpha$ is primitive, $t = 1, j = p$ and $s_0 = 1, s_1 = 2, ..., s_{p-1} = p$. The conclusion follows. $\square$

Let $G = \operatorname{PGL}_n(K)$ where $n = p^s \times t$, $\gcd(p, t) = 1$ and $s \geq 1$. For each factorization $n = p^r \times k$ where $1 \leq r \leq s$, $\bar{\Gamma}_r = \langle \bar{A}_0, \bar{B}_0, \bar{A}_1, \bar{B}_1, ..., \bar{A}_{r-1}, \bar{B}_{r-1} \rangle \leq G$. Let $D_r = \langle \bar{A}_0, \bar{A}_1, ..., \bar{A}_{r-1} \rangle$ and $B_r = \langle \bar{B}_0, \bar{B}_1, ..., \bar{B}_{r-1} \rangle$.

The next theorem gives the structure of the toral elementary abelian $p$-subgroups with disconnected centralizers.

**Theorem 4.27.** *Let $n = p^s \times t$ where $p$ is an odd prime, $\gcd(p, t) = 1$ and $s, t \geq 1$. The conjugacy classes of the toral elementary abelian $p$-subgroups of $\operatorname{PGL}_n(K)$ with disconnected centralizers have representatives $D_r \times H$ for $1 \leq r \leq s$, $n = p^r \times k$ and $H$ is either trivial or a representative of a conjugacy class of the toral elementary abelian $p$-subgroups of $\operatorname{PGL}_k(K)$ with $C_{\operatorname{PGL}_k(K)}(H) = C_{\operatorname{PGL}_k(K)}(H)^\circ$.*

*Proof.* The proof is similar to that for Theorem 4.9. □

**Theorem 4.28.** *Let $n = p^s \times t$ where $p$ is an odd prime, $\gcd(p, t) = 1$ and $s, t \geq 1$. For $n = p^r \times k$ where $1 \leq r \leq s$, if we let $D = D_r \times H$ in $G = \operatorname{PGL}_n(K)$ where $H$ is either trivial or a representative of a conjugacy class of the toral elementary abelian $p$-subgroups of $\operatorname{PGL}_k(K)$ with $C_{\operatorname{PGL}_k(K)}(H)$ connected, then*
$$C_G(D) = C_G(D)^\circ . B_r.$$

*Proof.* The proof is similar to that for Theorem 4.10. □

Now we study the normalizer quotient.

**Theorem 4.29.** *Let $n = p^s \times t$ where $p$ is an odd prime, $\gcd(p, t) = 1$ and $s, t \geq 1$. If $D \leqslant G = \operatorname{PGL}_n(K)$ is a toral elementary abelian $p$-subgroup with a disconnected centralizer, then*
$$N_G(D)/C_G(D) \geq \operatorname{GL}_r(p)$$
*where $1 \leq r \leq s$ and $n = p^r \times (p^{s-r} \times t)$. Hence there are two $F$-classes of $N_G(D)/C_G(D)^\circ$ in $C_G(D)/C_G(D)^\circ$, one consisting of the identity and the other the remaining elements.*

*Proof.* The proof is similar to that for Theorem 4.11. □

**Corollary 4.30.** *Let $G = \operatorname{PGL}_n(\mathbb{C})$ where $p^s \mid n$ but $p^{s+1} \nmid n$ for $p$ an odd prime and $s \geq 1$. Let $\bar{A}$ be a representative of a conjugacy class of the toral elementary abelian $p$-subgroups of rank $r \leq s$. If $C_G(\bar{A}) = C_G(\bar{A})^\circ$, then $\bar{A}$ has an element of order $p$ not conjugate to $e_{n/p}$.*

*Proof.* Suppose each nontrivial element of $\bar{A}$ is conjugate to $e_{n/p}$. In $\bar{\Gamma}_r$, we construct the toral elementary abelian $p$-subgroup $D_r = \langle \bar{A}_0, \bar{A}_1, ..., \bar{A}_{r-1} \rangle$. Since the conjugacy class containing $D_r$ is unique in $G$ and by Theorem 4.28, $C_G(D_r)$ is disconnected, we conclude that $\bar{A}$ and $D_r$ are conjugate in $G$ and $C_G(\bar{A})$ is disconnected, a contradiction. □

**Theorem 4.31.** *Let $G = \mathrm{PGL}_n(\mathbb{C})$ where $p \mid n$ and $p$ is an odd prime. Let $\bar{A}$ be a representative of a conjugacy class of the toral elementary abelian p-subgroups of $G$. If $C_G(\bar{A}) = C_G(\bar{A})^\circ$, then there exists a choice of generators of $\bar{A}$ such that no generator is conjugate to $e_{n/p}$.*

*Proof.* Let $n = p^u \times w$ where $\gcd(p, w) = 1$ and $u, w \geq 1$. First, not every nontrivial element of $\bar{A}$ is conjugate to $e_{n/p}$; if it holds, then $\bar{A}$ has rank at most $u$ and we apply Corollary 4.30 to reach a contradiction. Now suppose that every choice of generators $g_1, ..., g_s, g_{s+1}, ..., g_{s+t}$ of $\bar{A}$ satisfies the following: none of $g_1, g_2, ..., g_s$ is conjugate to $e_{n/p}$ but $xy$ is conjugate to $e_{n/p}$ for each $x \in H := \langle g_1, g_2, ..., g_s \rangle$ and for each $y \in \langle g_{s+1}, g_{s+2}, ..., g_{s+t} \rangle$.

We claim that if this holds, then $C_G(\bar{A})$ is disconnected. Let $E_2 = \langle \bar{A}_0, ..., \bar{A}_{t-1} \rangle \times H$ and $E_3 = \langle \bar{A}_0, \bar{B}_0, ..., \bar{A}_{t-1}, \bar{B}_{t-1} \rangle \times H$. Now $E_3$ is nontoral in $G$ and by Theorem 4.4,

$$N_G(E_3)/C_G(E_3) = \begin{pmatrix} \mathrm{Sp}_{2t}(p) & 0 \\ p^{s \times 2t} & W_{\mathrm{PGL}_k(C)}(H) \end{pmatrix}.$$

We conclude that $sx$ is conjugate to $e_{n/p}$ for $s \in \langle \bar{A}_0, ..., \bar{A}_{t-1} \rangle$ and for $x \in H$. Since the conjugacy class containing $\langle \bar{A}_0, ..., \bar{A}_{t-1} \rangle$ is unique in $G$, $\bar{A}$ and $E_2$ are conjugate in $G$ and by Theorem 4.28, $C_G(\bar{A})$ is not connected, a contradiction. $\square$

Now we obtain the precise structure of the normalizer quotient.

**Theorem 4.32.** *Let $G = \mathrm{PGL}_n(K)$ where $n = p^s \times t$, $p$ is an odd prime, $\gcd(p, t) = 1$ and $s, t \geq 1$. For $n = p^r \times (p^{s-r} \times t)$ where $1 \leq r \leq s$, if $D = D_r \times \bar{U}$ where $\bar{U} = \bar{U}_1 \otimes I_{p^r}$ and $\bar{U}_1$ is either trivial or a representative of a conjugacy class of toral elementary abelian p-subgroups of $\mathrm{PGL}_{p^{s-r}t}(K)$ with $C_{\mathrm{PGL}_{p^{s-r}t}(K)}(\bar{U}_1)$ connected, then*

$$N_G(D)/C_G(D) \cong p^{\mathrm{rk}\bar{U} \times r} : (\mathrm{GL}_r(p) \times W_{\mathrm{PGL}_{p^{s-r}t}(K)}(\bar{U}_1)).$$

*Proof.* The proof is similar to that for Theorem 4.14. $\square$

4.6.3. *Descending to the finite groups.*

Let $\Lambda$ be a nontoral elementary abelian p-subgroup of $\mathrm{PGL}_n(K)$ where $n = p^r \times k$ and both $r$ and $k$ are positive integers. Let $\mathrm{PGL}_n(K)^F = \mathrm{PGL}_n(q)$ where $q \equiv 1 \bmod p$ and $\Lambda \leq \mathrm{PGL}_n(K)^F$. We denote $N_{\mathrm{PGL}_n(K)}(\Lambda)$ by $N$ and $C_{\mathrm{PGL}_n(K)}(\Lambda)$ by $C$. Let $E$ be the nontoral elementary abelian p-subgroup of $\mathrm{PGL}_n(\mathbb{C})$ corresponding to $\Lambda$ under the bijection of Proposition 3.1. By Theorem 4.4, $E = \bar{\Gamma}_r \times \bar{A}$.

Now when $\bar{A}$ is trivial, $E = \bar{\Gamma}_r$, and $C_{\mathrm{PGL}_n(\mathbb{C})}(E) \cong \bar{\Gamma}_r \times \mathrm{PGL}_k(\mathbb{C})$ and $W_{\mathrm{PGL}_n(\mathbb{C})}(E) = \mathrm{Sp}_{2r}(p)$. By Theorem 4.23, for $r \geq 1$, the group $N/C^\circ \cong p^{2r}.\mathrm{Sp}_{2r}(p)$ is a subgroup of $\mathrm{PGL}_{p^r}(q)$ when $q \equiv 1 \bmod p$ and is therefore centralized by $F$. Hence, when $\Lambda = p^{2r}$, there are two $F$-classes of $N/C^\circ$ in $C/C^\circ$, one consisting of the identity and the other all the remaining elements. Correspondingly, there are two $\mathrm{PGL}_n(q)$-conjugacy classes of elementary abelian $p$-subgroups with representatives $E_1$ and $E_2$. Their local structure is the following:

$$\begin{aligned}
C_{\mathrm{PGL}_n(q)}(E_1) &\cong \bar{\Gamma}_r \times \mathrm{PGL}_k(q) \\
N_{\mathrm{PGL}_n(q)}(E_1) &\cong (\bar{\Gamma}_r \times \mathrm{PGL}_k(q)).\mathrm{Sp}_{2r}(p) \\
C_{\mathrm{PGL}_n(q)}(E_2) &\cong \bar{\Gamma}_r \times \mathrm{PGL}_k(q) \\
N_{\mathrm{PGL}_n(q)}(E_2) &\cong (\bar{\Gamma}_r \times \mathrm{PGL}_k(q)).(p^{2r-1} : \mathrm{Sp}_{2(r-1)}(p)).
\end{aligned}$$

Next, we focus on the nontrivial $\bar{A}$ and consider two cases.

**Case 1**. $\bar{A}$ is nontrivial toral in $\mathrm{PGL}_k(\mathbb{C})$ with $C_{\mathrm{PGL}_k(\mathbb{C})}(\bar{A})$ connected.

In this case,

$$C_{\mathrm{PGL}_n(\mathbb{C})}(E)/C_{\mathrm{PGL}_n(\mathbb{C})}(E)^\circ \cong \bar{\Gamma}_r$$

and

$$W_{\mathrm{PGL}_n(\mathbb{C})}(E) = p^{\mathrm{rk}\bar{A} \times 2r} : (\mathrm{Sp}_{2r}(p) \times W_{\mathrm{PGL}_k(\mathbb{C})}(\bar{A})).$$

Let $\Lambda = p^{2r} \times \Lambda_w$ correspond to $E$ where $\Lambda_w$ has a connected centralizer in $\mathrm{PGL}_k(K)$. Therefore, there are two $F$-classes of $N/C^\circ$ in $C/C^\circ$, one consisting of the identity and the other all the remaining elements of $p^{2r}$, the first component of $\Lambda$. Correspondingly, there are two $\mathrm{PGL}_n(q)$-conjugacy classes of elementary abelian $p$-subgroups with representatives $E_1$ and $E_2$. Their local structure is the following:

$$\begin{aligned}
C_{\mathrm{PGL}_n(q)}(E_1) &\cong \bar{\Gamma}_r \times (C_{\mathrm{PGL}_k(K)}(\Lambda_w))^F \\
C_{\mathrm{PGL}_n(q)}(E_2) &\cong \bar{\Gamma}_r \times (C_{\mathrm{PGL}_k(K)}(\Lambda_w))^F \\
N_{\mathrm{PGL}_n(q)}(E_1) &\cong (\bar{\Gamma}_r \times (C_{\mathrm{PGL}_k(K)}(\Lambda_w))^F).W_1 \\
N_{\mathrm{PGL}_n(q)}(E_2) &\cong (\bar{\Gamma}_r \times (C_{\mathrm{PGL}_k(K)}(\Lambda_w))^F).W_2
\end{aligned}$$

where

$$\begin{aligned}
W_1 &= p^{\mathrm{rk}\bar{A} \times 2r} : (\mathrm{Sp}_{2r}(p) \times W_{\mathrm{PGL}_k(\mathbb{C})}(\bar{A})) \\
W_2 &= p^{\mathrm{rk}\bar{A} \times 2r} : ((p^{2r-1} : \mathrm{Sp}_{2(r-1)}(p)) \times W_{\mathrm{PGL}_k(\mathbb{C})}(\bar{A})).
\end{aligned}$$

**Case 2**. $\bar{A}$ is nontrivial toral in $\mathrm{PGL}_k(\mathbb{C})$ with $C_{\mathrm{PGL}_k(\mathbb{C})}(\bar{A})$ disconnected.

Let $G = \mathrm{PGL}_n(K)$. Let $\Lambda = p^{2r} \times \Lambda_w \leqslant G$ correspond to $E$ where $\Lambda_w$ is a representative of a conjugacy class of toral elementary abelian $p$-subgroups of $\mathrm{PGL}_k(K)$ with $C_{\mathrm{PGL}_k(K)}(\Lambda_w)$ disconnected.

Let $k = p^{r_1} \times k_1$ where both $r_1$ and $k_1$ are positive integers. By Theorem 4.27, $\Lambda_w = D_{r_1} \times H$. Here $D_{r_1} = \langle \bar{A}_0^k, ..., \bar{A}_{r_1-1}^k \rangle$ where the superscript means elements in $\mathrm{PGL}_k(K)$ and $H$ is either trivial or a representative of a conjugacy class of the toral elementary abelian $p$-subgroups of $\mathrm{PGL}_{k_1}(K)$ with $C_{\mathrm{PGL}_{k_1}(K)}(H)$ connected. Let $B_{r_1} = \langle \bar{B}_0^k, ..., \bar{B}_{r_1-1}^k \rangle$. Now

$$C_G(\Lambda) \cong \bar{\Gamma}_r \times (C_{\mathrm{PGL}_k(K)}(\Lambda_w)^\circ . B_{r_1})$$

$$N_G(\Lambda)/C_G(\Lambda) \cong p^{\mathrm{rk}\Lambda_w \times 2r} : (\mathrm{Sp}_{2r}(p) \times W_{\mathrm{PGL}_k(K)}(\Lambda_w))$$

where

$$W_{\mathrm{PGL}_k(K)}(\Lambda_w) \cong p^{\mathrm{rk}H \times r_1} : (\mathrm{GL}_{r_1}(p) \times W_{\mathrm{PGL}_{k_1}(K)}(H)).$$

We denote by $R_1$ the group $p^{\mathrm{rk}\Lambda_w \times 2r}$ in $N_G(\Lambda)/C_G(\Lambda)$. The element $\bar{A}_0^k$ of $D_{r_1}$ is a product of some generators of the maximal toral elementary abelian $p$-subgroup of $\mathrm{PGL}_k(K)$. Corresponding to this product, we construct $y \in R_1$ in the same way as in the proof of Proposition 4.25 to get

$$y = \left[ \Delta\left( \begin{pmatrix} I_{p^r} & & & & \\ & \Delta(\square)_{p^{r-1}} & & & \\ & & \Delta(\square^2)_{p^{r-1}} & & \\ & & & \ddots & \\ & & & & \Delta(\square^{p-1})_{p^{r-1}} \end{pmatrix} \right)_{\frac{k}{p}} \right] C_G(\Lambda)$$

where

$$\square = \begin{pmatrix} 1 & & & & \\ & \alpha & & & \\ & & \alpha^2 & & \\ & & & \ddots & \\ & & & & \alpha^{p-1} \end{pmatrix}.$$

Since

$$\bar{B}_0^k = \left[ \Delta\left( \begin{pmatrix} & & 1 & \\ & 1 & & \\ & \iddots & & \\ 1 & & & 1 \end{pmatrix} \right)_{\frac{k}{p}} \right],$$

we deduce that

$$y(\bar{B}_0^k \otimes I_{p^r})y^{-1} = (\bar{B}_0^k \otimes I_{p^r})\bar{A}_0^{-1}.$$

By Theorem 4.29 and the action of $R_1$ on $\bar{\Gamma}_r$, we conclude that if $\Lambda = p^{2r} \times \Lambda_w$ where $\Lambda_w$ has a disconnected centralizer in $\mathrm{PGL}_k(K)$, then there are three $F$-classes of $N/C^\circ$ in $C/C^\circ$. One consists of the identity, another consists of the nontrivial elements of $p^{2r}$ and the third consists of all the remaining elements with representative $c$. Correspondingly, there are three $\mathrm{PGL}_n(q)$-conjugacy classes of elementary abelian $p$-subgroups with representatives $E_1$, $E_2$ and $E_3$. Their local structure is the following:

$$\begin{aligned}
C_{\mathrm{PGL}_n(q)}(E_1) &\cong \bar{\Gamma}_r \times (((C_{\mathrm{PGL}_k(K)}(\Lambda_w)^\circ)^F . B_{r_1}) \\
C_{\mathrm{PGL}_n(q)}(E_2) &\cong \bar{\Gamma}_r \times (((C_{\mathrm{PGL}_k(K)}(\Lambda_w)^\circ)^F . B_{r_1}) \\
C_{\mathrm{PGL}_n(q)}(E_3) &\cong \bar{\Gamma}_r \times (((C_{\mathrm{PGL}_k(K)}(\Lambda_w)^\circ)^{cF} . B_{r_1}) \\
N_{\mathrm{PGL}_n(q)}(E_1) &\cong C_{\mathrm{PGL}_n(q)}(E_1).(p^{\mathrm{rk}\Lambda_w \times 2r} : (\mathrm{Sp}_{2r}(p) \times W_{\mathrm{PGL}_k(K)}(\Lambda_w))) \\
N_{\mathrm{PGL}_n(q)}(E_2) &\cong C_{\mathrm{PGL}_n(q)}(E_2).(p^{\mathrm{rk}\Lambda_w \times 2r} : ((p^{2r-1} : \mathrm{Sp}_{2(r-1)}(p)) \times W_{\mathrm{PGL}_k(K)}(\Lambda_w))) \\
N_{\mathrm{PGL}_n(q)}(E_3) &\cong C_{\mathrm{PGL}_n(q)}(E_3).(p^{(\mathrm{rk}\Lambda_w - 1) \times 2r} : (\mathrm{Sp}_{2r}(p) \times W_3))
\end{aligned}$$

where
$$W_3 \cong p^{\text{rk}H \times r_1} : (\begin{pmatrix} 1 & * \\ 0 & \text{GL}_{r_1-1}(p) \end{pmatrix} \times W_{\text{PGL}_{k_1}(K)}(H)).$$

We have finished the investigation of the nontoral elementary abelian $p$-subgroups of $\text{PGL}_n(q)$ for odd $p$. Next, we carry out the parallel work for $\text{PGU}_n(q)$.

### 4.7. Nontoral elementary abelian $p$-subgroups of $\text{PGU}_n(q)$ where $p$ is odd

In this section, we examine the elementary abelian $p$-subgroups of $\text{PGU}_n(q)$ for odd $p$. We obtain the local structure of nontoral elementary abelian $p$-subgroups in Theorem 4.33 and then we descend to the finite groups and discuss separately the cases where $\bar{A}$ has a connected and disconnected centralizer.

Let $q$ be a power of a prime $\ell$ where $p \neq \ell$. Let $K$ be an algebraically closed field of characteristic $\ell$ and let $F$ be a Steinberg endomorphism of $\text{PGL}_n(K)$ such that $\text{PGL}_n(K)^F = \text{PGU}_n(q)$.

By [3, 1B], we conclude that if $p^m \mid n$ and $q \equiv -1 \mod p$, then $p^{1+2m}.\text{Sp}_{2m}(p) \leq \text{GU}_n(q)$. The following result is used to study the $F$-action.

**Theorem 4.33.** *Let $p$ be an odd prime and $m \geq 1$. If there exists a nontoral elementary abelian subgroup $E$ of order $p^{2m}$ of $\text{PGU}_n(q)$ with*
$$C_{\text{PGU}_n(q)}(E) \cong E \times \text{PGU}_k(q)$$
*where $1 \leq k < n$ and*
$$N_{\text{PGU}_n(q)}(E) \cong C_{\text{PGU}_n(q)}(E).\text{Sp}_{2m}(p),$$
*then $p^m \mid n$ and $q \equiv -1 \mod p$. Conversely, if $p^m \mid n$ and $q \equiv -1 \mod p$, then there exists a nontoral elementary abelian subgroup $E$ of order $p^{2m}$ of $\text{PGU}_n(q)$ with*
$$C_{\text{PGU}_n(q)}(E) \cong E \times \text{PGU}_k(q)$$
*where $1 \leq k < n$ and*
$$N_{\text{PGU}_n(q)}(E) \cong C_{\text{PGU}_n(q)}(E).\text{Sp}_{2m}(p).$$

*Proof.* The proof is similar to that for Theorem 4.23. $\square$

Now we investigate the local structure and descend to the finite groups. Let $\Lambda$ be a nontoral elementary abelian $p$-subgroup of $\text{PGL}_n(K)$ where $n = p^r \times k$ where both $r$ and $k$ are positive integers,. We assume that $\text{PGL}_n(K)^F = \text{PGU}_n(q)$ where $q \equiv -1 \mod p$ and $\Lambda \leq \text{PGL}_n(K)^F$. We denote $N_{\text{PGL}_n(K)}(\Lambda)$ by $N$ and $C_{\text{PGL}_n(K)}(\Lambda)$ by $C$. By Theorem 4.33, $N/C^\circ$ is centralized by $F$ when $q \equiv -1 \mod p$.

Let $E$ be the nontoral elementary abelian $p$-subgroup of $\mathrm{PGL}_n(\mathbb{C})$ corresponding to $\Lambda$ under the bijection of Proposition 3.1. By Theorem 4.4,
$$E = \bar{\Gamma}_r \times \bar{A}$$
where $\bar{A}$ is either trivial or a representative of a conjugacy class of the toral elementary abelian $p$-subgroups of $\mathrm{PGL}_k(\mathbb{C})$. By Proposition 3.2, $\Lambda = p^{2r} \times \Lambda_w$ where $\Lambda_w$ is either trivial or a toral elementary abelian $p$-subgroup of $\mathrm{PGL}_k(K)$.

When $\bar{A}$ is trivial, $E = \bar{\Gamma}_r$, and $C_{\mathrm{PGL}_n(\mathbb{C})}(E) \cong \bar{\Gamma}_r \times \mathrm{PGL}_k(\mathbb{C})$ and $W_{\mathrm{PGL}_n(\mathbb{C})}(E) = \mathrm{Sp}_{2r}(p)$. By Theorem 4.33, for $r \geq 1$, the group $\bar{\Gamma}_r.\mathrm{Sp}_{2r}(p)$ is a subgroup of $\mathrm{PGU}_{p^r}(q)$ for $q \equiv -1 \bmod p$ and is therefore centralized by $F$. Hence, when $\Lambda = p^{2r}$, there are two $F$-classes of $N/C^\circ$ in $C/C^\circ$, one consisting of the identity and the other all the remaining elements. Correspondingly, there are two $\mathrm{PGU}_n(q)$-conjugacy classes of elementary abelian $p$-subgroups with representatives $E_1$ and $E_2$. Their local structure is the following:

$$\begin{aligned}
C_{\mathrm{PGU}_n(q)}(E_1) &\cong \bar{\Gamma}_r \times \mathrm{PGU}_k(q) \\
N_{\mathrm{PGU}_n(q)}(E_1) &\cong (\bar{\Gamma}_r \times \mathrm{PGU}_k(q)).\mathrm{Sp}_{2r}(p) \\
C_{\mathrm{PGU}_n(q)}(E_2) &\cong \bar{\Gamma}_r \times \mathrm{PGU}_k(q) \\
N_{\mathrm{PGU}_n(q)}(E_2) &\cong (\bar{\Gamma}_r \times \mathrm{PGU}_k(q)).(p^{2r-1} : \mathrm{Sp}_{2(r-1)}(p)).
\end{aligned}$$

Next, we focus on the nontrivial $\bar{A}$ and consider two cases.

**Case 1**. $\bar{A}$ is nontrivial toral in $\mathrm{PGL}_k(\mathbb{C})$ with $C_{\mathrm{PGL}_k(\mathbb{C})}(\bar{A})$ connected.

In this case,
$$C_{\mathrm{PGL}_n(\mathbb{C})}(E)/C_{\mathrm{PGL}_n(\mathbb{C})}(E)^\circ \cong \bar{\Gamma}_r$$
and
$$W_{\mathrm{PGL}_n(\mathbb{C})}(E) = p^{\mathrm{rk}\bar{A} \times 2r} : (\mathrm{Sp}_{2r}(p) \times W_{\mathrm{PGL}_k(\mathbb{C})}(\bar{A})).$$

Let $\Lambda = p^{2r} \times \Lambda_w$ correspond to $E$ where $\Lambda_w$ has a connected centralizer in $\mathrm{PGL}_k(K)$. Therefore, there are two $F$-classes of $N/C^\circ$ in $C/C^\circ$, one consisting of the identity and the other all the remaining elements. Correspondingly, there are two $\mathrm{PGU}_n(q)$-conjugacy classes of elementary abelian $p$-subgroups with representatives $E_1$ and $E_2$. Their local structure is the following:

$$\begin{aligned}
C_{\mathrm{PGU}_n(q)}(E_1) &\cong \bar{\Gamma}_r \times (C_{\mathrm{PGL}_k(K)}(\Lambda_w))^F \\
C_{\mathrm{PGU}_n(q)}(E_2) &\cong \bar{\Gamma}_r \times (C_{\mathrm{PGL}_k(K)}(\Lambda_w))^F \\
N_{\mathrm{PGU}_n(q)}(E_1) &\cong (\bar{\Gamma}_r \times (C_{\mathrm{PGL}_k(K)}(\Lambda_w))^F).W_1 \\
N_{\mathrm{PGU}_n(q)}(E_2) &\cong (\bar{\Gamma}_r \times (C_{\mathrm{PGL}_k(K)}(\Lambda_w))^F).W_2
\end{aligned}$$

where
$$W_1 = p^{\mathrm{rk}\bar{A}\times 2r} : (\mathrm{Sp}_{2r}(p) \times W_{\mathrm{PGL}_k(\mathbb{C})}(\bar{A}))$$
$$W_2 = p^{\mathrm{rk}\bar{A}\times 2r} : ((p^{2r-1} : \mathrm{Sp}_{2(r-1)}(p)) \times W_{\mathrm{PGL}_k(\mathbb{C})}(\bar{A})).$$

**Case 2.** $\bar{A}$ is toral in $\mathrm{PGL}_k(\mathbb{C})$ with $C_{\mathrm{PGL}_k(\mathbb{C})}(\bar{A})$ disconnected.

Let $G = \mathrm{PGL}_n(K)$. Let $\Lambda = p^{2r} \times \Lambda_w \leqslant G$ correspond to $E$ where $\Lambda_w$ is a representative of a conjugacy class of toral elementary abelian $p$-subgroups of $\mathrm{PGL}_k(K)$ with $C_{\mathrm{PGL}_k(K)}(\Lambda_w)$ disconnected.

Let $k = p^{r_1} \times k_1$ where both $r_1$ and $k_1$ are positive integers. By Theorem 4.27, $\Lambda_w = D_{r_1} \times H$. Here $D_{r_1} = \langle \bar{A}_0^k, ..., \bar{A}_{r_1-1}^k \rangle$ where the superscript means elements in $\mathrm{PGL}_k(K)$ and $H$ is either trivial or a representative of a conjugacy class of the toral elementary abelian $p$-subgroups of $\mathrm{PGL}_{k_1}(K)$ with $C_{\mathrm{PGL}_{k_1}(K)}(H)$ connected. Let $B_{r_1} = \langle \bar{B}_0^k, ..., \bar{B}_{r_1-1}^k \rangle$. Now

$$C_G(\Lambda) \cong \bar{\Gamma}_r \times (C_{\mathrm{PGL}_k(K)}(\Lambda_w)^\circ . B_{r_1})$$

$$N_G(\Lambda)/C_G(\Lambda) \cong p^{\mathrm{rk}\Lambda_w \times 2r} : (\mathrm{Sp}_{2r}(p) \times W_{\mathrm{PGL}_k(K)}(\Lambda_w))$$

where

$$W_{\mathrm{PGL}_k(K)}(\Lambda_w) \cong p^{\mathrm{rk}H \times r_1} : (\mathrm{GL}_{r_1}(p) \times W_{\mathrm{PGL}_{k_1}(K)}(H)).$$

There are three $F$-classes of $N/C^\circ$ in $C/C^\circ$. One consists of the identity, another consists of the nontrivial elements of $p^{2r}$ and the third consists of all the remaining elements with representative $c$. Correspondingly, there are three $\mathrm{PGU}_n(q)$-conjugacy classes of elementary abelian $p$-subgroups with representatives $E_1$, $E_2$ and $E_3$. Their local structure is the following:

$$
\begin{aligned}
C_{\mathrm{PGU}_n(q)}(E_1) &\cong \bar{\Gamma}_r \times ((C_{\mathrm{PGL}_k(K)}(\Lambda_w)^\circ)^F . B_{r_1}) \\
C_{\mathrm{PGU}_n(q)}(E_2) &\cong \bar{\Gamma}_r \times ((C_{\mathrm{PGL}_k(K)}(\Lambda_w)^\circ)^F . B_{r_1}) \\
C_{\mathrm{PGU}_n(q)}(E_3) &\cong \bar{\Gamma}_r \times ((C_{\mathrm{PGL}_k(K)}(\Lambda_w)^\circ)^{cF} . B_{r_1}) \\
N_{\mathrm{PGU}_n(q)}(E_1) &\cong C_{\mathrm{PGU}_n(q)}(E_1) . (p^{\mathrm{rk}\Lambda_w \times 2r} : (\mathrm{Sp}_{2r}(p) \times W_{\mathrm{PGL}_k(K)}(\Lambda_w))) \\
N_{\mathrm{PGU}_n(q)}(E_2) &\cong C_{\mathrm{PGU}_n(q)}(E_2) . (p^{\mathrm{rk}\Lambda_w \times 2r} : ((p^{2r-1} : \mathrm{Sp}_{2(r-1)}(p)) \times W_{\mathrm{PGL}_k(K)}(\Lambda_w))) \\
N_{\mathrm{PGU}_n(q)}(E_3) &\cong C_{\mathrm{PGU}_n(q)}(E_3) . (p^{(\mathrm{rk}\Lambda_w - 1) \times 2r} : (\mathrm{Sp}_{2r}(p) \times W_3))
\end{aligned}
$$

where

$$W_3 \cong p^{\mathrm{rk}H \times r_1} : (\begin{pmatrix} 1 & * \\ 0 & \mathrm{GL}_{r_1-1}(p) \end{pmatrix} \times W_{\mathrm{PGL}_{k_1}(K)}(H)).$$

## Acknowledgments

This work forms part of my PhD [12] obtained from the University of Auckland under the supervision of Professors Jianbei An and Eamonn O'Brien. I am grateful for their assistance, and thank the University of Auckland for its financial support.

(Fu) DEPARTMENT OF MATHEMATICS, UNIVERSITY OF AUCKLAND, NZ

*Email address*: gerrardmzh@yahoo.com